\documentclass[12pt,reqno]{amsart}

\newtheorem{Theorem}{Theorem}
\newtheorem{Proposition}[Theorem]{Proposition}
\newtheorem{Corollary}[Theorem]{Corollary}
\newtheorem{Lemma}[Theorem]{Lemma}
\numberwithin{Theorem}{section}
\numberwithin{equation}{section}

\allowdisplaybreaks

\begin{document}

\title[A new $_1\psi_1$ summation in $A_n$]
{A new $\hbox{$\boldsymbol A_{\boldsymbol n}$}$ extension of
Ramanujan's $\hbox{${}_{\boldsymbol 1}\boldsymbol\psi_{\boldsymbol 1}$}$
summation with applications to multilateral
$\hbox{$\boldsymbol A_{\boldsymbol n}$}$ series
}
\author[Stephen C. Milne]{Stephen C. Milne$^*$}
\author{Michael Schlosser}

\address{Department of Mathematics, The Ohio State University,
231 West 18th Avenue, Columbus, Ohio 43210, USA}
\email{milne@math.ohio-state.edu, mschloss@math.ohio-state.edu}
\date{October 11, 2000}
\thanks{$^*$ S.~C.~Milne was partially supported by 
National Security Agency grant MDA904--99--1--0003.}

\dedicatory{{\small\rm Department of Mathematics,
The Ohio State University,\\
231 West 18th Avenue, Columbus, Ohio 43210, USA\\
{\rm E-mail:}} {\footnotesize\tt milne@math.ohio-state.edu,
mschloss@math.ohio-state.edu}}

\begin{abstract}
In this article, we derive some identities for multilateral basic
hypergeometric series associated to the root system $A_n$.
First, we apply Ismail's~\cite{ismail} argument to
an $A_n$ $q$-binomial theorem of Milne~\cite[Th.~5.42]{milne}
and derive a new $A_n$ generalization of Ramanujan's
$_1\psi_1$ summation theorem. From this new $A_n$ $_1\psi_1$ summation
and from an $A_n$ $_1\psi_1$ summation of Gustafson~\cite{gusmult}
we deduce two lemmas for deriving simple $A_n$ generalizations of
bilateral basic hypergeometric series identities. These lemmas are
closely related to the Macdonald identities for $A_n$.
As samples for possible applications of these lemmas, we provide several
$A_n$ extensions of Bailey's $_2\psi_2$ transformations,
and several $A_n$ extensions of a particular $_2\psi_2$ summation.
\end{abstract}

\maketitle

\section{Introduction}

The theory of basic hypergeometric series (cf.~\cite{grhyp})
consists of many known summation and transformation formulas.
The most important of these is probably the $q$-binomial
theorem, a summation first discovered by Cauchy~\cite{cauchy}.
Surpringly, the $q$-binomial theorem admits a bilateral
generalization, the $_1\psi_1$ summation theorem, first discovered
by Ramanujan~\cite{hardy}.
Other important identities for basic hypergeometric series
include the $q$-Gau\ss{} summation and Heine's $_2\phi_1$ transformations.
These and many other basic hypergeometric series identities conspicuously
appear in combinatorics and in
related areas such as number theory, statistics, physics,
and representation theory of Lie algebras, see Andrews~\cite{qandrews}.

Multiple basic hypergeometric series
associated to the root system $A_n$ (or equivalently, associated
to the unitary group $U\!(n+1)$) have been investigated
by various authors.
Many different types of such series exist
in the literature. The multivariable series we consider in this article
have their origin in the work of the three mathematical physicists
Biedenharn, Holman and Louck, see \cite{holman} and \cite{holbl}.
Their work was done in the context of the quantum theory of angular
momentum, using methods relying on the representation theory of $U\!(n)$.
In the sequel, substantial developments have taken place.
Extensive investigations in the theory of
multiple basic hypergeometric series associated to the root system
$A_n$ have been carried out by Gustafson, Milne, and their co-workers.
As result, many of the classical formulas for basic hypergeometric series
(cf.~\cite{grhyp}) have already been generalized to the setting of $A_n$
series. For some selected results on multiple basic hypergeometric series
associated to $A_n$, see the references
\cite{bhatmil}, \cite{dengus}, \cite{gusmult}, \cite{gus}, \cite{milmac},
\cite{milram}, \cite{milbas}, \cite{milq}, \cite{miln1}, \cite{milwhip},
\cite{milne}, \cite{millil2}, \cite{milnew}, \cite{schlossmmi},
\cite{schlnammi}, and \cite{schlnmmi}.

There are different methods for obtaining identities for
$A_n$ basic hypergeometric series. Partial fraction decompositions
and $q$-difference equations are often involved in initially
deriving such identities (see, e.g., \cite{bhatmil}, \cite{gus},
and \cite{milmac}). Further, where summations for multidimensional basic
hypergeometric series are already known, multidimensional matrix
inversions can often be utilized for obtaining new summation theorems
for multidimensional basic hypergeometric series
(see~\cite{bhatmil}, \cite{milne}, \cite{millil2},
\cite{schlossmmi}, \cite{schlnammi}, and \cite{schlnmmi}).
But there is also another, simpler, way of obtaining identities
for $A_n$ basic hypergeometric series. By utilizing Lemma~7.3 of
Milne~\cite{milne}, see Lemma~\ref{lemmilne} in this article,
and by using identities of the classical
one-dimensional theory, simple identities for $A_n$ series
can be derived.

In this article, we find two multilateral generalizations
of \cite[Lemma~7.3]{milne}, see Lemmas~\ref{lem1} and \ref{lem2}.
These lemmas are closely related to the Macdonald~\cite{macd}
identities for the affine root system $A_n$.
By using our lemmas combined with bilateral one-dimensional series
identities we are able to derive simple multilateral identities for
$A_n$ series. We give some particular applications of this method.
The $A_n$ $_2\psi_2$ transformations and summations
given in this article are just
samples of the possible applications.
It must be said that the identities obtained by this method 
concern $A_n$ series of ``simpler type" and are apparently not as deep
as many of the $A_n$ identities in the above mentioned references.
Nevertheless, in spite of, or maybe even because of the ``simplicity"
of these $A_n$ series our formulas could be useful in future applications. 

Our article is organized as follows:
In Section~\ref{secproc}, we introduce some notation and give some
background information.
In Section~\ref{secqbin}, we apply Ismail's~\cite{ismail}
analytic continuation argument
to an $A_n$ $q$-binomial theorem of Milne~\cite[Theorem~5.42]{milne}
to derive a new $A_n$ extension of Ramanujan's~\cite{hardy} $_1\psi_1$
summation theorem. In \cite{milram} a similar argument
was used to find the first $U\!(n)$ generalization of the $_1\psi_1$
summation. More recently, motivated by \cite{milschur},
Kaneko~\cite{kaneko} utilized this type of argument to derive a
$_1\psi_1$ summation theorem for multiple basic hypergeometric series
of Macdonald polynomial argument.
In Section~\ref{seclemma}, we deduce from our new $A_n$ $_1\psi_1$
summation and from Gustafson's~\cite[Theorem~1.17]{gusmult}
$_1\psi_1$ summation two lemmas for deriving simple
multilateral series identities in $A_n$. We discuss the connection
of these lemmas with the Macdonald identities for $A_n$,
partly following the similar analysis of \cite{milram}. 
Finally, in Section~\ref{secbail}, we apply these lemmas to classical
(one-dimensional) formulas for $_2\psi_2$ series. As result, we
deduce several (different) $A_n$ extensions of
Bailey's~\cite{bail22} $_2\psi_2$ transformations, and moreover,
deduce several (different) $A_n$ extensions of a particular
summation for $_2\psi_2$ series.

\section{Background and notation}\label{secproc}

Let us first recall some standard basic hypergeometric notation 
(cf.~\cite{grhyp}). Let $q$ be a complex number such that $0<|q|<1$.
We define the {\em $q$-shifted factorial} for all integers $k$ by
\begin{equation*}
(a)_{\infty}\equiv(a;q)_{\infty}:=\prod_{j=0}^{\infty}(1-aq^j)
\qquad\text{and}\qquad
(a)_k\equiv(a;q)_k:=\frac{(a)_{\infty}}{(aq^k)_{\infty}}.
\end{equation*}
For brevity, we employ the usual notation
\begin{equation*}
(a_1,\ldots,a_m)_k\equiv (a_1)_k\dots(a_m)_k
\end{equation*}
where $k$ is an integer or infinity. Further, we utilize the notations
\begin{equation}\label{defhyp}
_r\phi_s\!\left[\begin{matrix}a_1,a_2,\dots,a_r\\
b_1,b_2,\dots,b_s\end{matrix};q,z\right]:=
\sum _{k=0} ^{\infty}\frac {(a_1,a_2,\dots,a_r)_k}
{(q,b_1,\dots,b_s)_k}\left((-1)^kq^{\binom k2}\right)^{1+s-r}z^k,
\end{equation}
and
\begin{equation}\label{defhypb}
_r\psi_s\!\left[\begin{matrix}a_1,a_2,\dots,a_r\\
b_1,b_2,\dots,b_s\end{matrix};q,z\right]:=
\sum _{k=-\infty} ^{\infty}\frac {(a_1,a_2,\dots,a_r)_k}
{(b_1,b_2,\dots,b_s)_k}\left((-1)^kq^{\binom k2}\right)^{s-r}z^k,
\end{equation}
for {\em basic hypergeometric $_r\phi_s$ series}, and {\em bilateral basic
hypergeometric $_r\psi_s$
series}, respectively. See \cite[p.~25 and p.~125]{grhyp} for the criteria
of when these series terminate, or, if not, when they converge. 
In this article, we make use of some of the elementary
identities for $q$-shifted factorials, listed in \cite[Appendix I]{grhyp}.

Next, we note the convention for naming the multiple series in this
article as $A_n$ basic hypergeometric series.
We consider multiple series of the form 
\begin{equation*}
\sum_{k_1,\dots,k_n=-\infty}^{\infty} S({\mathbf k}),
\end{equation*}
where ${\mathbf k}=(k_1,\dots,k_n)$,
which reduce to classical basic hypergeometric series when $n=1$.
Such a multiple series 
is called an $A_n$ basic hypergeometric series if the summand
$S({\mathbf k})$ contains the factor
\begin{equation}\label{anvandy}
\prod_{1\le i<j\le n} \left(\frac {x_iq^{k_i}-x_jq^{k_j}}
{x_i-x_j}\right).
\end{equation} 
A typical example is the left hand side of \eqref{an1phi0gl}.
A reason for naming these series as $A_n$ series is that \eqref{anvandy}
is closely associated with the product side of the Weyl denominator
formula for the root system $A_n$, see \cite{bhat} and \cite{stan}.

For multidimensional series, we also employ the notation
$|{\mathbf k}|$ for $(k_1+\dots+k_n)$ where ${\mathbf k}=(k_1,\dots,k_n)$.
The convergence of multiple series
can be checked by application of the multiple power series ratio
test~\cite{horn}.
For explicit examples of how to use the multiple power series ratio
test, see \cite[Sec.~5]{milne}.

\section{An $A_n$ extension of Ramanujan's $_1\psi_1$ summation}
\label{secqbin}

One of the most important summation theorems for basic hypergeometric series
is the classical $q$-binomial theorem (cf.~\cite[Eq.~(II.3)]{grhyp}),
\begin{equation}\label{10gl}
{}_1\phi_0\!\left[\begin{matrix}a\\
-\end{matrix};q,z\right]=
\frac{(az)_{\infty}}{(z)_{\infty}},
\end{equation}
where $|z|<1$.

A bilateral extension of \eqref{10gl} is Ramanujan's~\cite{hardy}
$_1\psi_1$ summation theorem (cf.~\cite[Eq.~(5.2.1)]{grhyp}), 
\begin{equation}\label{11gl}
{}_1\psi_1\!\left[\begin{matrix}a\\
b\end{matrix};q,z\right]=
\frac{(q,b/a,az,q/az)_{\infty}}{(b,q/a,z,b/az)_{\infty}},
\end{equation}
where $|b/a|<|z|<1$. Clearly, the $b=q$ case of \eqref{11gl} is \eqref{10gl}.

Theorem~5.42 of \cite{milne} is one of the many multivariable
generalizations of \eqref{10gl}. It can be stated as follows:

\begin{Theorem}[An $A_n$ $q$-binomial theorem]\label{an1phi0}
Let $a$, $x_1,\dots,x_n$, and $z$ be indeterminate,
let $n\ge 1$, and suppose that
none of the denominators in \eqref{an1phi0gl} vanishes. Then
\begin{multline}\label{an1phi0gl}
\sum_{k_1,\dots,k_n=0}^{\infty}\Bigg(
\prod_{1\le i<j\le n}\left(\frac {x_iq^{k_i}-x_jq^{k_j}}{x_i-x_j}\right)
\prod_{i,j=1}^n\left(\frac{x_i}{x_j}q\right)_{k_i}^{-1}
\prod_{i=1}^nx_i^{nk_i-|{\mathbf k}|}\\\times
(a)_{|{\mathbf k}|}(-1)^{(n-1)|{\mathbf k}|}
q^{-\binom{|{\mathbf k}|}2+n\sum_{i=1}^n\binom{k_i}2}
z^{|\mathbf k|}\Bigg)
=\frac{(az)_{\infty}}{(z)_{\infty}},
\end{multline}
provided $|z|<\left|q^{\frac
{n-1}2}x_j^{-n}\prod_{i=1}^nx_i\right|$
for $j=1,\dots,n$.
\end{Theorem}

We now apply Ismail's~\cite{ismail} argument and extend
Theorem~\ref{an1phi0} to

\begin{Theorem}[An $A_n$ $_1\psi_1$ summation]\label{an1psi1}
Let $a$, $b_1,\dots,b_n$, $x_1,\dots,x_n$, and $z$ be indeterminate,
let $n\ge 1$, and suppose that
none of the denominators in \eqref{an1psi1gl} vanishes. Then
\begin{multline}\label{an1psi1gl}
\sum_{k_1,\dots,k_n=-\infty}^{\infty}\Bigg(
\prod_{1\le i<j\le n}\left(\frac {x_iq^{k_i}-x_jq^{k_j}}{x_i-x_j}\right)
\prod_{i,j=1}^n\left(\frac{x_i}{x_j}b_j\right)_{k_i}^{-1}
\prod_{i=1}^nx_i^{nk_i-|{\mathbf k}|}\\\times
(a)_{|{\mathbf k}|}(-1)^{(n-1)|{\mathbf k}|}
q^{-\binom{|{\mathbf k}|}2+n\sum_{i=1}^n\binom{k_i}2}
z^{|\mathbf k|}\Bigg)\\
=\frac{(az,q/az,b_1\dots b_nq^{1-n}/a)_{\infty}}
{(z,b_1\dots b_nq^{1-n}/az,q/a)_{\infty}}
\prod_{i,j=1}^n\frac{\left(\frac{x_i}{x_j}q\right)_{\infty}}
{\left(\frac{x_i}{x_j}b_j\right)_{\infty}},
\end{multline}
provided $|b_1\dots b_nq^{1-n}/a|<|z|<\left|q^{\frac
{n-1}2}x_j^{-n}\prod_{i=1}^nx_i\right|$
for $j=1,\dots,n$.
\end{Theorem}
\begin{proof}
We apply Ismail's argument successively to the parameters $b_1,\dots,b_n$
using \eqref{an1phi0gl}. The multiple series identity in \eqref{an1psi1gl}
is analytic in each of the parameters $b_1,\dots,b_n$ in a domain around
the origin. Now, the identity is true for $b_1=q^{1+m_1},
b_2=q^{1+m_2},\dots,$ and $b_n=q^{1+m_n}$, by the $A_n$ $q$-binomial theorem
in Theorem~\ref{an1phi0} (see below for the details). This holds for
all $m_1,\dots,m_n\ge 0$.
Since $\lim_{m_1\to\infty}q^{1+m_1}=0$ is an interior point
in the domain of analycity of $b_1$, by analytic continuation,
we obtain an identity for $b_1$.
By iterating this argument for $b_2,\dots,b_n$, we establish
\eqref{an1psi1gl} for general $b_1,\dots,b_n$.

The details are displayed as follows. Setting $b_i=q^{1+m_i}$, for
$i=1,\dots,n$, the left side of \eqref{an1psi1gl} becomes
\begin{multline}\label{longl1}
\sum_{\begin{smallmatrix}-m_i\le k_i\le\infty\\
i=1,\dots,n\end{smallmatrix}}\Bigg(
\prod_{1\le i<j\le n}\left(\frac {x_iq^{k_i}-x_jq^{k_j}}{x_i-x_j}\right)
\prod_{i,j=1}^n\left(\frac{x_i}{x_j}q^{1+m_j}\right)_{k_i}^{-1}
\prod_{i=1}^nx_i^{nk_i-|{\mathbf k}|}\\\times
(a)_{|{\mathbf k}|}(-1)^{(n-1)|{\mathbf k}|}
q^{-\binom{|{\mathbf k}|}2+n\sum_{i=1}^n\binom{k_i}2}
z^{|\mathbf k|}\Bigg).
\end{multline}
We shift the summation indices in \eqref{longl1} by
$k_i\mapsto k_i-m_i$, for $i=1,\dots,n$ and obtain
\begin{multline}\label{longl}
q^{-\binom{|{\mathbf m}|+1}2+n\sum_{i=1} ^n\binom{m_i+1}2}
(-1)^{(n-1)|{\mathbf m}|}(a)_{-|\mathbf m|}z^{-|\mathbf m|}
\prod_{i=1}^nx_i^{|{\mathbf m}|-nm_i}
\prod_{i,j=1}^n\left(\frac{x_i}{x_j}q^{1+m_j}\right)_{-m_i}^{-1}\\
\times\sum_{k_1,\dots,k_n=0}^{\infty}\Bigg(
\prod_{1\le i<j\le n}\left(\frac {x_iq^{-m_i+k_i}-x_jq^{-m_j+k_j}}
{x_i-x_j}\right)
\prod_{i,j=1}^n\left(\frac{x_i}{x_j}q^{1+m_j-m_i}\right)_{k_i}^{-1}
\\\times
(aq^{-|{\mathbf m}|})_{|{\mathbf k}|}(-1)^{(n-1)|{\mathbf k}|}
q^{-\binom{|{\mathbf k}|}2+n\sum_{i=1}^n\binom{k_i}2}
z^{|\mathbf k|}
\prod_{i=1}^n\left(x_iq^{-m_i}\right)^{nk_i-|{\mathbf k}|}
\Bigg)\\
=q^{n\sum_{i=1} ^n\binom{m_i+1}2}(-1)^{n|{\mathbf m}|}
(az)^{-|\mathbf m|}(q/a)_{|\mathbf m|}^{-1}\\\times
\prod_{i=1}^nx_i^{|{\mathbf m}|-nm_i}
\prod_{i,j=1}^n\frac{\left(\frac{x_i}{x_j}q\right)_{m_j}}
{\left(\frac{x_i}{x_j}q\right)_{m_j-m_i}}
\prod_{1\le i<j\le n}\left(\frac {x_iq^{-m_i}-x_jq^{-m_j}}
{x_i-x_j}\right)\\\times
\sum _{k_1,\dots,k_n=0}^{\infty}\Bigg(
\prod_{1\le i<j\le n}\left(\frac {x_iq^{-m_i+k_i}-x_jq^{-m_j+k_j}}
{x_iq^{-m_i}-x_jq^{-m_j}}\right)
\prod_{i,j=1}^n\left(\frac{x_i}{x_j}q^{1+m_j-m_i}\right)_{k_i}^{-1}
\\\times
(aq^{-|{\mathbf m}|})_{|{\mathbf k}|}(-1)^{(n-1)|{\mathbf k}|}
q^{-\binom{|{\mathbf k}|}2+n\sum_{i=1}^n\binom{k_i}2}z^{|\mathbf k|}
\prod_{i=1}^n\left(x_iq^{-m_i}\right)^{nk_i-|{\mathbf k}|}\Bigg).
\end{multline}
Next, we apply the $y_i\mapsto-m_i$, $i=1,\dots,n$, case of
\cite[Lemma~3.12]{milne}, specifically
\begin{multline}
\prod_{i,j=1}^n\left(\frac{x_i}{x_j}q\right)_{m_j-m_i}=
(-1)^{(n-1)|{\mathbf m}|}
q^{-\binom{|\mathbf m|+1}2+n\sum_{i=1} ^n\binom{m_i+1}2}\\\times
\prod_{i=1}^nx_i^{|{\mathbf m}|-nm_i}
\prod_{1\le i<j\le n}\left(\frac {x_iq^{-m_i}-x_jq^{-m_j}}
{x_i-x_j}\right),
\end{multline}
and the $a\mapsto aq^{-|\mathbf m|}$,
$x_i\mapsto x_iq^{-m_i}$, $i=1,\dots,n$, case of the multidimensional
summation theorem in \eqref{an1phi0gl} to simplify the expression obtained
in \eqref{longl} to
\begin{equation*}
q^{\binom{|\mathbf m|+1}2}(-az)^{-|\mathbf m|}
\frac{(azq^{-|\mathbf m|})_{\infty}}
{(q/a)_{|\mathbf m|}(z)_{\infty}}
\prod_{i,j=1}^n\left(\frac{x_i}{x_j}q\right)_{m_j}.
\end{equation*}
Now, this can easily be further transformed into
\begin{equation*}
\frac{(q^{1+|\mathbf m|}/a,az,q/az)_{\infty}}
{(q/a,z,q^{1+|\mathbf m|}/az)_{\infty}}
\prod_{i,j=1}^n\frac{\left(\frac{x_i}{x_j}q\right)_{\infty}}
{\left(\frac{x_i}{x_j}q^{1+m_j}\right)_{\infty}},
\end{equation*}
which is exactly the $b_i=q^{1+m_i}$, $i=1,\dots,n$, case of the right side
of \eqref{an1psi1gl}.
\end{proof}

If we set $z\mapsto -z/a$, and $b_i=0$, $i=1,\dots,n$ in \eqref{an1psi1gl},
and then let $a\to\infty$, we obtain an $A_n$ generalization of Jacobi's
triple product identity, equivalent to Theorem~3.7 of \cite{milram}.

\section{Two lemmas for deriving multilateral $A_n$ series identities}
\label{seclemma}

As an immediate consequence of a fundamental theorem for $A_n$
series~\cite[Theorem~1.49]{milmac}, the first author~\cite[Lemma~7.3]{milne}
of this article derived the following lemma, which is
\begin{Lemma}[Milne]\label{lemmilne}
Let $a_1,\dots,a_n$ and $x_1,\dots,x_n$ be indeterminate, let $N$ be
a nonnegative integer, let $n\ge 1$, and suppose that
none of the denominators in \eqref{lemmilnegl} vanishes.
Then, if $f(m)$ is an arbitrary function of nonnegative integers $m$, we have
\begin{multline}\label{lemmilnegl}
\sum_{m=0}^N\frac{(a_1a_2\dots a_n)_m}{(q)_m}f(m)\\
=\sum_{\begin{smallmatrix}k_1,\dots,k_n\ge 0\\
0\le |{\mathbf k}|\le N\end{smallmatrix}}
\prod_{1\le i<j\le n}\left(\frac {x_iq^{k_i}-x_jq^{k_j}}{x_i-x_j}\right)
\prod_{i,j=1}^n\frac{(\frac{x_i}{x_j}a_j)_{k_i}}{(\frac{x_i}{x_j}q)_{k_i}}
\cdot f(|{\mathbf k}|).
\end{multline}
\end{Lemma}
With Lemma~\ref{lemmilne} and one-dimensional basic hypergeometric
series identities, (simple) identities for $A_n$ series can be derived.
Some examples are given in \cite[Sec.~7]{milne}.

In this section, we provide two new lemmas, see Lemmas~\ref{lem1} and
\ref{lem2}, which similarly can be used for deriving simple $A_n$
generalizations of {\em bilateral} basic hypergeometric series identities.
We make use of our $A_n$ extension of Ramanujan's $_1\psi_1$ summation
in Theorem~\ref{an1psi1} and of an $A_n$ $_1\psi_1$ summation
by Gustafson~\cite[Theorem~1.17]{gusmult}, see Theorem~\ref{an1psi1g}.

Since, for $|b_1\dots b_nq^{1-n}/a|<|z|<1$,
\begin{equation*}
{}_1\psi_1\!\left[\begin{matrix}a\\
b_1\dots b_nq^{1-n}\end{matrix};q,z\right]=
\frac{(q,b_1\dots b_nq^{1-n}/a,az,q/az)_{\infty}}
{(b_1\dots b_nq^{1-n},q/a,z,b_1\dots b_nq^{1-n}/az)_{\infty}},
\end{equation*}
by Ramanujan's $_1\psi_1$ summation \eqref{11gl},
we immediately see from \eqref{an1psi1gl} that
\begin{multline}\label{an1psi1gla}
\sum_{k_1,\dots,k_n=-\infty}^{\infty}\Bigg(
\prod_{1\le i<j\le n}\left(\frac {x_iq^{k_i}-x_jq^{k_j}}{x_i-x_j}\right)
\prod_{i,j=1}^n\left(\frac{x_i}{x_j}b_j\right)_{k_i}^{-1}
\prod_{i=1}^nx_i^{nk_i-|{\mathbf k}|}\\\times
(a)_{|{\mathbf k}|}(-1)^{(n-1)|{\mathbf k}|}
q^{-\binom{|{\mathbf k}|}2+n\sum_{i=1}^n\binom{k_i}2}
z^{|\mathbf k|}\Bigg)\\
=\frac{(b_1\dots b_nq^{1-n})_{\infty}}{(q)_{\infty}}
\prod_{i,j=1}^n\frac{\left(\frac{x_i}{x_j}q\right)_{\infty}}
{\left(\frac{x_i}{x_j}b_j\right)_{\infty}}
\sum_{k=-\infty}^{\infty}\frac{(a)_k}{(b_1\dots b_nq^{1-n})_k}z^k
\end{multline}
(provided $|z|<1$ and $|b_1\dots b_nq^{1-n}/a|<|z|<\left|q^{\frac
{n-1}2}x_j^{-n}\prod_{i=1}^nx_i\right|$
for $j=1,\dots,n$).

In \eqref{an1psi1gla}, we equate coefficients of
$(a)_mz^m$ and extract
\begin{Proposition}\label{an1psi1N}
Let $b_1,\dots,b_n$ and $x_1,\dots,x_n$ be indeterminate,
let $m$ be an integer, let $n\ge 1$, and suppose that none of the
denominators in \eqref{an1psi1glN} vanishes. Then
\begin{multline}\label{an1psi1glN}
\sum_{\begin{smallmatrix}-\infty\le k_1,\dots,k_n\le\infty\\
|{\mathbf k}|=m\end{smallmatrix}}\Bigg(
\prod_{1\le i<j\le n}\left(\frac {x_iq^{k_i}-x_jq^{k_j}}{x_i-x_j}\right)
\prod_{i,j=1}^n\left(\frac{x_i}{x_j}b_j\right)_{k_i}^{-1}
\prod_{i=1}^nx_i^{nk_i-|{\mathbf k}|}\\\times
(-1)^{(n-1)|{\mathbf k}|}
q^{-\binom{|{\mathbf k}|}2+n\sum_{i=1}^n\binom{k_i}2}\Bigg)\\
=\frac{(b_1\dots b_nq^{1-n})_{\infty}}{(q)_{\infty}}
\prod_{i,j=1}^n\frac{\left(\frac{x_i}{x_j}q\right)_{\infty}}
{\left(\frac{x_i}{x_j}b_j\right)_{\infty}}\,\cdot
\frac{1}{(b_1\dots b_nq^{1-n})_m}.
\end{multline}
\end{Proposition}
We state Proposition~\ref{an1psi1N} although it
is just a special case of Proposition~\ref{an1psi1gN}.
We utilize the $m=0$ case of Proposition~\ref{an1psi1N}
in the proof of Theorem~\ref{an2psi2s}.

Now, if we multiply both sides of \eqref{an1psi1glN} by
\begin{equation*}
\frac{(q)_{\infty}}{(b_1\dots b_nq^{1-n})_{\infty}}
\prod_{i,j=1}^n\frac{\left(\frac{x_i}{x_j}b_j\right)_{\infty}}
{\left(\frac{x_i}{x_j}q\right)_{\infty}}\cdot f(m),
\end{equation*}
for suitable $f(m)$, and sum over all integers $m$, we obtain
\begin{Lemma}\label{lem1}
Let $b_1,\dots,b_n$, $x_1,\dots,x_n$ be indeterminate, let $n\ge 1$,
and suppose that none of the denominators in \eqref{lem1gl} vanishes.
Then, if $f(m)$ is an arbitrary function of integers $m$, we have
\begin{multline}\label{lem1gl}
\sum_{m=-\infty}^{\infty}\frac{f(m)}{(b_1\dots b_nq^{1-n})_m}
=\frac{(q)_{\infty}}{(b_1\dots b_nq^{1-n})_{\infty}}
\prod_{i,j=1}^n\frac{\left(\frac{x_i}{x_j}b_j\right)_{\infty}}
{\left(\frac{x_i}{x_j}q\right)_{\infty}}\\\times
\sum_{k_1,\dots,k_n=-\infty}^{\infty}
\Bigg(\prod_{1\le i<j\le n}\left(\frac {x_iq^{k_i}-x_jq^{k_j}}{x_i-x_j}\right)
\prod_{i,j=1}^n\left(\frac{x_i}{x_j}b_j\right)_{k_i}^{-1}
\prod_{i=1}^nx_i^{nk_i-|{\mathbf k}|}\\\times
(-1)^{(n-1)|{\mathbf k}|}
q^{-\binom{|{\mathbf k}|}2+n\sum_{i=1}^n\binom{k_i}2}
\cdot f(|{\mathbf k}|)\Bigg),
\end{multline}
provided the series converge.
\end{Lemma}

Thus, with Lemma~\ref{lem1}, we can use one-dimensional bilateral
series identities to obtain identities for multilateral $A_n$ series.

The special case of Lemma~\ref{lem1}, where $b_i=q$, for $i=1,\dots,n$
is worth noting:
\begin{Corollary}\label{cor1}
Let $x_1,\dots,x_n$ be indeterminate, let $n\ge 1$,
and suppose that none of the denominators in \eqref{cor1gl} vanishes.
Then, if $f(m)$ is an arbitrary function of nonnegative integers $m$, we have
\begin{multline}\label{cor1gl}
\sum_{m=0}^{\infty}\frac{f(m)}{(q)_m}
=\sum_{k_1,\dots,k_n=0}^{\infty}
\Bigg(\prod_{1\le i<j\le n}\left(\frac {x_iq^{k_i}-x_jq^{k_j}}{x_i-x_j}\right)
\prod_{i,j=1}^n\left(\frac{x_i}{x_j}q\right)_{k_i}^{-1}
\prod_{i=1}^nx_i^{nk_i-|{\mathbf k}|}\\\times
(-1)^{(n-1)|{\mathbf k}|}
q^{-\binom{|{\mathbf k}|}2+n\sum_{i=1}^n\binom{k_i}2}
\cdot f(|{\mathbf k}|)\Bigg),
\end{multline}
provided the series converge.
\end{Corollary}

Corollary~\ref{cor1} can be also obtained by specializing
Lemma~\ref{lemmilne}. Namely, by setting
\begin{equation*}
f(m)\mapsto (-1)^mq^{-\binom m2}(a_1a_2\dots a_n)^{-m}f(m)
\end{equation*}
in Lemma~\ref{lemmilne},
and then letting $N\to\infty$ and $a_i\to\infty$, for $i=1,\dots,n$,
we also obtain Corollary~\ref{cor1}. 

Next, we put our attention towards the derivation of another
lemma for deriving multilateral series identities.
For this, we utilize Gustafson's~\cite[Theorem~1.17]{gusmult}
multivariable generalization of Ramanujan's $_1\psi_1$ summation~\eqref{11gl}.
\begin{Theorem}[(Gustafson) An $A_n$ $_1\psi_1$ summation]\label{an1psi1g}
Let $a_1,\dots,a_n$, $b_1,\dots,b_n$, $z$, and $x_1,\dots,x_n$ be
indeterminate, let $n\ge 1$, and suppose that
none of the denominators in \eqref{an1psi1ggl} vanishes. Then
\begin{multline}\label{an1psi1ggl}
\sum_{k_1,\dots,k_n=-\infty}^{\infty}
\prod_{1\le i<j\le n}\left(\frac {x_iq^{k_i}-x_jq^{k_j}}{x_i-x_j}\right)
\prod_{i,j=1}^n\frac{\left(\frac{x_i}{x_j}a_j\right)_{k_i}}
{\left(\frac{x_i}{x_j}b_j\right)_{k_i}}z^{|\mathbf k|}\\
=\frac{(a_1\dots a_nz,q/a_1\dots a_nz)_{\infty}}
{(z,b_1\dots b_nq^{1-n}/a_1\dots a_nz)_{\infty}}
\prod_{i,j=1}^n
\frac{\left(\frac{x_i}{x_j}q,\frac{x_ib_j}{x_ja_i}\right)_{\infty}}
{\left(\frac{x_i}{x_j}b_j,\frac{x_iq}{x_ja_i}\right)_{\infty}},
\end{multline}
provided $|b_1\dots b_nq^{1-n}/a_1\dots a_n|<|z|<1$.
\end{Theorem}

Since, for $|b_1\dots b_nq^{1-n}/a_1\dots a_n|<|z|<1$,
\begin{equation*}
{}_1\psi_1\!\left[\begin{matrix}a_1\dots a_n\\
b_1\dots b_nq^{1-n}\end{matrix};q,z\right]=
\frac{(q,b_1\dots b_nq^{1-n}/a_1\dots a_n,a_1\dots a_nz,
q/a_1\dots a_nz)_{\infty}}
{(b_1\dots b_nq^{1-n},q/a_1\dots a_n,z,
b_1\dots b_nq^{1-n}/a_1\dots a_nz)_{\infty}},
\end{equation*}
by Ramanujan's $_1\psi_1$ summation \eqref{11gl},
we immediately see from \eqref{an1psi1ggl} that
\begin{multline}\label{an1psi1ggla}
\sum_{k_1,\dots,k_n=-\infty}^{\infty}
\prod_{1\le i<j\le n}\left(\frac {x_iq^{k_i}-x_jq^{k_j}}{x_i-x_j}\right)
\prod_{i,j=1}^n\frac{\left(\frac{x_i}{x_j}a_j\right)_{k_i}}
{\left(\frac{x_i}{x_j}b_j\right)_{k_i}}z^{|\mathbf k|}\\
=\frac{(b_1\dots b_nq^{1-n},q/a_1\dots a_n)_{\infty}}
{(q,b_1\dots b_nq^{1-n}/a_1\dots a_n)_{\infty}}
\prod_{i,j=1}^n
\frac{\left(\frac{x_i}{x_j}q,\frac{x_ib_j}{x_ja_i}\right)_{\infty}}
{\left(\frac{x_i}{x_j}b_j,\frac{x_iq}{x_ja_i}\right)_{\infty}}
\sum_{k=-\infty}^{\infty}\frac{(a_1\dots a_n)_k}{(b_1\dots b_nq^{1-n})_k}z^k
\end{multline}
(provided $|b_1\dots b_nq^{1-n}/a_1\dots a_n|<|z|<1$).

In \eqref{an1psi1ggla}, we equate coefficients of
$z^m$ and extract
\begin{Proposition}\label{an1psi1gN}
Let $a_1,\dots,a_n$, $b_1,\dots,b_n$, and $x_1,\dots,x_n$ be indeterminate,
let $m$ be an integer, let $n\ge 1$, and suppose that none of the
denominators in \eqref{an1psi1gglN} vanishes. Then
\begin{multline}\label{an1psi1gglN}
\sum_{\begin{smallmatrix}-\infty\le k_1,\dots,k_n\le\infty\\
|{\mathbf k}|=m\end{smallmatrix}}
\prod_{1\le i<j\le n}\left(\frac {x_iq^{k_i}-x_jq^{k_j}}{x_i-x_j}\right)
\prod_{i,j=1}^n\frac{\left(\frac{x_i}{x_j}a_j\right)_{k_i}}
{\left(\frac{x_i}{x_j}b_j\right)_{k_i}}\\
=\frac{(b_1\dots b_nq^{1-n},q/a_1\dots a_n)_{\infty}}
{(q,b_1\dots b_nq^{1-n}/a_1\dots a_n)_{\infty}}
\prod_{i,j=1}^n
\frac{\left(\frac{x_i}{x_j}q,\frac{x_ib_j}{x_ja_i}\right)_{\infty}}
{\left(\frac{x_i}{x_j}b_j,\frac{x_iq}{x_ja_i}\right)_{\infty}}\,
\cdot\frac{(a_1\dots a_n)_m}{(b_1\dots b_nq^{1-n})_m},
\end{multline}
provided $|b_1\dots b_nq^{1-n}/a_1\dots a_n|<1$.
\end{Proposition}
The $b_i=b$, $i=1,\dots,n$, case of Proposition~\ref{an1psi1gN} was
established in \cite[Theorem~1.21]{milram}.

A specialization of Proposition~\ref{an1psi1gN} gives
Proposition~\ref{an1psi1N}.
Namely, if we divide both sides of \eqref{an1psi1gglN} by
$(a_1\dots a_n)_m$ and then let $a_i\to\infty$, $i=1,\dots,n$, we
obtain \eqref{an1psi1glN}.

The $m=0$ case of Proposition~\ref{an1psi1gN} was established by
Gustafson in \cite[Theorem~1.15]{gusmult}:
\begin{Theorem}[(Gustafson) An $A_{n-1}$ $_6\psi_6$
summation]\label{an1psi1cN}
Let $a_1,\dots,a_n$, $b_1,\dots,b_n$, and $x_1,\dots,x_n$ be indeterminate,
let $n\ge 1$, and suppose that none of the
denominators in \eqref{an1psi1cglN} vanishes. Then
\begin{multline}\label{an1psi1cglN}
\sum_{\begin{smallmatrix}-\infty\le k_1,\dots,k_n\le\infty\\
|{\mathbf k}|=0\end{smallmatrix}}
\prod_{1\le i<j\le n}\left(\frac {x_iq^{k_i}-x_jq^{k_j}}{x_i-x_j}\right)
\prod_{i,j=1}^n\frac{\left(\frac{x_i}{x_j}a_j\right)_{k_i}}
{\left(\frac{x_i}{x_j}b_j\right)_{k_i}}\\
=\frac{(b_1\dots b_nq^{1-n},q/a_1\dots a_n)_{\infty}}
{(q,b_1\dots b_nq^{1-n}/a_1\dots a_n)_{\infty}}
\prod_{i,j=1}^n
\frac{\left(\frac{x_i}{x_j}q,\frac{x_ib_j}{x_ja_i}\right)_{\infty}}
{\left(\frac{x_i}{x_j}b_j,\frac{x_iq}{x_ja_i}\right)_{\infty}},
\end{multline}
provided $|b_1\dots b_nq^{1-n}/a_1\dots a_n|<1$.
\end{Theorem}
The $n=2$ case of Theorem~\ref{an1psi1cN} is equivalent to
Bailey's~\cite[Eq.~(4.7)]{bail66} very-well-poised $_6\psi_6$ summation
(cf.~\cite[Eq.~(5.3.1)]{grhyp}).

We utilize Theorem~\ref{an1psi1cN} in the proof of Theorem~\ref{an2psi23s}.

From Theorem~\ref{an1psi1cN} we immediately deduce
a $_1\psi_1$/$_6\psi_6$ generalization of
the Macdonald identities for $A_n$, generalizing
Theorem~1.24 of \cite{milram}.
The analysis is similar to that in \cite{milram}
where the $b_i=b$, $i=1,\dots,n$, case of Theorem~\ref{an1psi1cN}
was utilized to obtain \cite[Theorem~1.24]{milram}.
The following result appears implicitly in \cite[Sec.~7]{gus}. 

\begin{Theorem}[(Gustafson) A $_1\psi_1$ generalization of the Macdonald
identities for $A_n$]\label{an1psi1mN}
Let $a_1,\dots,a_n$, $b_1,\dots,b_n$, and $x_1,\dots,x_n$ be indeterminate,
let $n\ge 1$, and suppose that none of the
denominators in \eqref{an1psi1mglN} vanishes. Then
\begin{multline}\label{an1psi1mglN}
\sum_{\sigma\in{\mathcal S}_n}\varepsilon(\sigma)
\prod_{i=1}^nx_{\sigma(i)}^{i-\sigma(i)}
\sum_{\begin{smallmatrix}-\infty\le k_1,\dots,k_n\le\infty\\
|{\mathbf k}|=0\end{smallmatrix}}
q^{\sum_{i=1}^r(i-1)k_{\sigma(i)}}
\prod_{i,j=1}^n\frac{\left(\frac{x_i}{x_j}a_j\right)_{k_i}}
{\left(\frac{x_i}{x_j}b_j\right)_{k_i}}\\
=\frac{(b_1\dots b_nq^{1-n},q/a_1\dots a_n)_{\infty}}
{(q,b_1\dots b_nq^{1-n}/a_1\dots a_n)_{\infty}}
\prod_{1\le i<j\le n}\left(1-\frac{x_i}{x_j}\right)
\prod_{i,j=1}^n
\frac{\left(\frac{x_i}{x_j}q,\frac{x_ib_j}{x_ja_i}\right)_{\infty}}
{\left(\frac{x_i}{x_j}b_j,\frac{x_iq}{x_ja_i}\right)_{\infty}},
\end{multline}
provided $|b_1\dots b_nq^{1-n}/a_1\dots a_n|<1$,
where ${\mathcal S}_n$ is the symmetric group of order $n$,
and $\varepsilon(\sigma)$ is the sign of the permutation $\sigma$. 
\end{Theorem}

Replacing $a_i$ and $b_i$ by $-1/c$ and $0$, respectively, for
$i=1,\dots,n$, in Theorem~\ref{an1psi1mN}, simplifying and then
letting $c\to 0$ yields Equation (4.3) of \cite{milmac} which
is equivalent to the Macdonald identities for $A_n$, see
\cite[Sec.~4]{milmac}. Thus, Theorem~\ref{an1psi1mN} may be viewed as
a generalization of the Macdonald identities for $A_n$ with the extra
parameters $a_1,\dots,a_n$, and $b_1,\dots,b_n$.

For future reference, we write down the $b_i=a_iq$, $i=1,\dots,n$,
case of Theorems~\ref{an1psi1mN} and \ref{an1psi1gN}. Note that
this case is valid since the convergence condition
$|b_1\dots b_nq^{1-n}/a_1\dots a_n|<1$ becomes $|q|<1$.
After a routine simplification, \eqref{an1psi1mglN} becomes
\begin{multline}\label{an1psi1mnglN}
\sum_{\sigma\in{\mathcal S}_n}\varepsilon(\sigma)
\prod_{i=1}^nx_{\sigma(i)}^{i-\sigma(i)}
\sum_{\begin{smallmatrix}-\infty\le k_1,\dots,k_n\le\infty\\
|{\mathbf k}|=0\end{smallmatrix}}
q^{\sum_{i=1}^r(i-1)k_{\sigma(i)}}
\prod_{i,j=1}^n\frac{\left(1-\frac{x_i}{x_j}a_j\right)}
{\left(1-\frac{x_i}{x_j}a_jq^{k_i}\right)}\\
=\frac{(a_1\dots a_nq,q/a_1\dots a_n)_{\infty}}
{(q,q)_{\infty}}
\prod_{1\le i<j\le n}\left(1-\frac{x_i}{x_j}\right)
\prod_{i,j=1}^n
\frac{\left(\frac{x_i}{x_j}q,\frac{x_ia_j}{x_ja_i}q\right)_{\infty}}
{\left(\frac{x_i}{x_j}a_jq,\frac{x_iq}{x_ja_i}\right)_{\infty}}.
\end{multline}
Similarly, \eqref{an1psi1gglN} becomes
\begin{multline}\label{an1psi1gnglN}
\sum_{\begin{smallmatrix}-\infty\le k_1,\dots,k_n\le\infty\\
|{\mathbf k}|=m\end{smallmatrix}}
\prod_{1\le i<j\le n}\left(\frac {x_iq^{k_i}-x_jq^{k_j}}{x_i-x_j}\right)
\prod_{i,j=1}^n\frac{\left(1-\frac{x_i}{x_j}a_j\right)}
{\left(1-\frac{x_i}{x_j}a_jq^{k_i}\right)}\\
=\frac{(a_1\dots a_n,q/a_1\dots a_n)_{\infty}}
{(1-a_1\dots a_nq^m)\,(q,q)_{\infty}}
\prod_{i,j=1}^n
\frac{\left(\frac{x_i}{x_j}q,\frac{x_ia_j}{x_ja_i}q\right)_{\infty}}
{\left(\frac{x_i}{x_j}a_jq,\frac{x_iq}{x_ja_i}\right)_{\infty}}.
\end{multline}
Equations \eqref{an1psi1mnglN} and \eqref{an1psi1gnglN} extend
(3.16) and (3.17) of \cite{milram}, respectively, to which they reduce when
$a_i=a$, for $i=1,\dots,n$.

Now, let us return to our objective of finding a multilateral generalization
of Lemma~\ref{lemmilne}. If we multiply both sides of \eqref{an1psi1gglN} by
\begin{equation*}
\frac{(q,b_1\dots b_nq^{1-n}/a_1\dots a_n)_{\infty}}
{(b_1\dots b_nq^{1-n},q/a_1\dots a_n)_{\infty}}
\prod_{i,j=1}^n\frac{\left(\frac{x_i}{x_j}b_j,
\frac{x_iq}{x_ja_i}\right)_{\infty}}
{\left(\frac{x_i}{x_j}q,\frac{x_ib_j}{x_ja_i}\right)_{\infty}}\cdot g(m),
\end{equation*}
for suitable $g(m)$, and sum over all integers $m$, we obtain
\begin{Lemma}\label{lem2}
Let $a_1,\dots,a_n$, $b_1,\dots,b_n$, $x_1,\dots,x_n$ be indeterminate,
let $n\ge 1$, and suppose that none of the denominators in
\eqref{lem2gl} vanishes.
Then, if $g(m)$ is an arbitrary function of integers $m$, we have
\begin{multline}\label{lem2gl}
\sum_{m=-\infty}^{\infty}\frac{(a_1\dots a_n)_m}{(b_1\dots b_nq^{1-n})_m}
g(m)\\
=\frac{(q,b_1\dots b_nq^{1-n}/a_1\dots a_n)_{\infty}}
{(b_1\dots b_nq^{1-n},q/a_1\dots a_n)_{\infty}}
\prod_{i,j=1}^n\frac{\left(\frac{x_i}{x_j}b_j,
\frac{x_iq}{x_ja_i}\right)_{\infty}}
{\left(\frac{x_i}{x_j}q,\frac{x_ib_j}{x_ja_i}\right)_{\infty}}\\\times
\sum_{k_1,\dots,k_n=-\infty}^{\infty}
\prod_{1\le i<j\le n}\left(\frac {x_iq^{k_i}-x_jq^{k_j}}{x_i-x_j}\right)
\prod_{i,j=1}^n\frac{\left(\frac{x_i}{x_j}a_j\right)_{k_i}}
{\left(\frac{x_i}{x_j}b_j\right)_{k_i}}
\cdot g(|{\mathbf k}|),
\end{multline}
provided the series converge.
\end{Lemma}

Hence, besides Lemma~\ref{lem1}, we can also use Lemma~\ref{lem2}
with one-dimensional bilateral
series identities to obtain identities for multilateral $A_n$ series.
Lemma~\ref{lem2} generalizes the $N\to\infty$ case of
Lemma~\ref{lemmilne} by additional parameters $b_1,\dots,b_n$, since
the special case $b_i=q$, for $i=1,\dots,n$, of Lemma~\ref{lem2}
boils down to the $N\to\infty$ case of Lemma~\ref{lemmilne}.

\section{Applications: Some $_2\psi_2$ formulas in $A_n$}
\label{secbail}

In this section, we illustrate the usefulness of the lemmas of
the preceding section and provide some multidimensional extensions
of Bailey's~\cite{bail22} $_2\psi_2$ transformations,
associated to the root system $A_n$. Further, as interesting special
cases of these $_2\psi_2$ transformations in $A_n$, we provide some
$_2\psi_2$ summations in $A_n$.

Using Ramanujan's $_1\psi_1$ summation \eqref{11gl} and elementary
manipulations of series, Bailey~\cite[Eq.~(2.3)]{bail22} derived the
transformation
\begin{equation}\label{22tgl}
{}_2\psi_2\!\left[\begin{matrix}a,b\\
c,d\end{matrix};q,z\right]
=\frac{(az,d/a,c/b,dq/abz)_{\infty}}{(z,d,q/b,cd/abz)_{\infty}}\:
{}_2\psi_2\!\left[\begin{matrix}a,abz/d\\
az,c\end{matrix};q,\frac da\right],
\end{equation}
where $\max(|z|,|cd/abz|,|d/a|,|c/b|)<1$.

Bailey's $_2\psi_2$ transformation can be iterated. The result is
\cite[Eq.~(2.4)]{bail22}
\begin{equation}\label{22tgl1}
{}_2\psi_2\!\left[\begin{matrix}a,b\\
c,d\end{matrix};q,z\right]
=\frac{(az,bz,cq/abz,dq/abz)_{\infty}}{(q/a,q/b,c,d)_{\infty}}\:
{}_2\psi_2\!\left[\begin{matrix}abz/c,abz/d\\
az,bz\end{matrix};q,\frac{cd}{abz}\right],
\end{equation}
where $\max(|z|,|cd/abz|)<1$.

We can specialize \eqref{22tgl} (or \eqref{22tgl1}) to obtain
a summation theorem for a particular $_2\psi_2$ series.
If $d=bq$ and $z=q/a$ in \eqref{22tgl}, then the series on the right side
reduces just to one term, $1$, and we have the summation
\begin{equation}\label{22sgl}
{}_2\psi_2\!\left[\begin{matrix}a,b\\
c,bq\end{matrix};q,\frac qa\right]
=\frac{(q,q,bq/a,c/b)_{\infty}}{(q/a,bq,q/b,c)_{\infty}},
\end{equation}
where $\max(|q/a|,|c|)<1$.

In the following subsections, we combine our
Lemmas~\ref{lem1} and \ref{lem2} from Section~\ref{seclemma}
together with the above one-dimensional $_2\psi_2$ formulas.
In Subsection~\ref{subsec1}, we derive several
multivariable extensions of Bailey's $_2\psi_2$ transformation formulas
\eqref{22tgl} and \eqref{22tgl1}.
In Subsection~\ref{subsec2} we derive multivariable extensions of
the $_2\psi_2$ summation in \eqref{22sgl}.

\subsection{Some $\hbox{$\boldsymbol A_{\boldsymbol n}$}$ extensions of
Bailey's $\hbox{${}_{\boldsymbol 2}\boldsymbol\psi_{\boldsymbol 2}$}$
transformations}
\label{subsec1}

We give several (but not all) of the possible $A_n$ $_2\psi_2$ transformations
which arise from Lemmas~\ref{lem1} and \ref{lem2}.

We start with two multivariable extensions of \eqref{22tgl}
which arise from Lemma~\ref{lem1}.

\begin{Theorem}[An $A_n$ $_2\psi_2$ transformation]\label{an2psi2}
Let $a$, $b$, $c_1,\dots,c_n$, $d$, $x_1,\dots,x_n$, $y_1,\dots,y_n$,
and $z$ be indeterminate, let $n\ge 1$, and suppose that
none of the denominators in \eqref{an2psi2gl} vanishes. Then
\begin{multline}\label{an2psi2gl}
\sum_{k_1,\dots,k_n=-\infty}^{\infty}\Bigg(
\prod_{1\le i<j\le n}\left(\frac {x_iq^{k_i}-x_jq^{k_j}}{x_i-x_j}\right)
\prod_{i,j=1}^n\left(\frac{x_i}{x_j}c_j\right)_{k_i}^{-1}
\prod_{i=1}^nx_i^{nk_i-|{\mathbf k}|}\\\times
\frac{(a,b)_{|{\mathbf k}|}}{(d)_{|{\mathbf k}|}}(-1)^{(n-1)|{\mathbf k}|}
q^{-\binom{|{\mathbf k}|}2+n\sum_{i=1}^n\binom{k_i}2}
z^{|\mathbf k|}\Bigg)\\
=\frac{(az,d/a,c_1\dots c_nq^{1-n}/b,dq/abz)_{\infty}}
{(z,d,q/b,c_1\dots c_ndq^{1-n}/abz)_{\infty}}
\prod_{i,j=1}^n
\frac{\left(\frac{x_i}{x_j}q,\frac{y_i}{y_j}c_j\right)_\infty}
{\left(\frac{y_i}{y_j}q,\frac{x_i}{x_j}c_j\right)_\infty}\\\times
\sum_{k_1,\dots,k_n=-\infty}^{\infty}\Bigg(
\prod_{1\le i<j\le n}\left(\frac {y_iq^{k_i}-y_jq^{k_j}}{y_i-y_j}\right)
\prod_{i,j=1}^n\left(\frac{y_i}{y_j}c_j\right)_{k_i}^{-1}
\prod_{i=1}^ny_i^{nk_i-|{\mathbf k}|}\\\times
\frac{(a,abz/d)_{|{\mathbf k}|}}{(az)_{|{\mathbf k}|}}(-1)^{(n-1)|{\mathbf k}|}
q^{-\binom{|{\mathbf k}|}2+n\sum_{i=1}^n\binom{k_i}2}
\left(\frac da\right)^{|\mathbf k|}\Bigg),
\end{multline}
provided $|c_1\dots c_ndq^{1-n}/ab|<|z|<\left|q^{\frac
{n-1}2}x_j^{-n}\prod_{i=1}^nx_i\right|$ and
$|c_1\dots c_ndq^{1-n}/ab|<|d/a|<\left|q^{\frac
{n-1}2}y_j^{-n}\prod_{i=1}^ny_i\right|$
for $j=1,\dots,n$.
\end{Theorem}
\begin{proof}
We have, for
$\max(|z|,|c_1\dots c_ndq^{1-n}/abz|,|d/a|,|c_1\dots c_nq^{1-n}/b|)<1$,
\begin{multline}\label{22t1gl}
{}_2\psi_2\!\left[\begin{matrix}a,b\\
c_1\dots c_nq^{1-n},d\end{matrix};q,z\right]\\
=\frac{(az,d/a,c_1\dots c_nq^{1-n}/b,dq/abz)_{\infty}}
{(z,d,q/b,c_1\dots c_ndq^{1-n}/abz)_{\infty}}\:
{}_2\psi_2\!\left[\begin{matrix}a,abz/d\\
az,c_1\dots c_nq^{1-n}\end{matrix};q,\frac da\right],
\end{multline}
by Bailey's $_2\psi_2$ transformation in \eqref{22tgl}.
Now we apply Lemma~\ref{lem1} to the $_2\psi_2$'s on the left and on the
right side of this transformation.
Specifically, we rewrite the $_2\psi_2$ on left side of \eqref{22t1gl}
by the $b_i\mapsto c_i$, $i=1,\dots,n$, and
\begin{equation*}
f(m)=\frac{(a,b)_m}{(d)_m}z^m
\end{equation*}
case of Lemma~\ref{lem1}.
The $_2\psi_2$ on the right side of \eqref{22t1gl} is rewritten
by the $b_i\mapsto c_i$, $x_i\mapsto y_i$, $i=1,\dots,n$, and
\begin{equation*}
f(m)=\frac{(a,abz/d)_m}{(az)_m}\left(\frac da\right)^m
\end{equation*}
case of Lemma~\ref{lem1}.
Finally, we divide both sides of the resulting equation by
\begin{equation}\label{divgl}
\frac{(q)_{\infty}}{(c_1\dots c_nq^{1-n})_{\infty}}
\prod_{i,j=1}^n\frac{\left(\frac{x_i}{x_j}c_j\right)_{\infty}}
{\left(\frac{x_i}{x_j}q\right)_{\infty}}
\end{equation}
and simplify to obtain \eqref{an2psi2gl}.
\end{proof}

\begin{Theorem}[An $A_n$ ${}_2\psi_2$ transformation]\label{an2psi22}
Let $a_1,a_2,\dots,a_n$, $b$, $c_1,\dots,c_n$, $d$, $x_1,\dots,x_n$,
$y_1,\dots,y_n$, and $z_1,\dots,z_n$ be indeterminate, let $n\ge 1$,
and suppose that none of the denominators in \eqref{an2psi22gl}
vanishes. Write $A\equiv a_1\dots a_n$, $C\equiv c_1\dots c_n$,
and $Z\equiv z_1\dots z_n$, for short. Then
\begin{multline}\label{an2psi22gl}
\sum_{k_1,\dots,k_n=-\infty}^{\infty}\Bigg(
\prod_{1\le i<j\le n}\left(\frac {x_iq^{k_i}-x_jq^{k_j}}{x_i-x_j}\right)
\prod_{i,j=1}^n\left(\frac{x_i}{x_j}c_j\right)_{k_i}^{-1}
\prod_{i=1}^nx_i^{nk_i-|{\mathbf k}|}\\\times
\frac{(Aq^{1-n},b)_{|{\mathbf k}|}}
{(d)_{|{\mathbf k}|}}(-1)^{(n-1)|{\mathbf k}|}
q^{-\binom{|{\mathbf k}|}2+n\sum_{i=1}^n\binom{k_i}2}
Z^{|\mathbf k|}\Bigg)\\
=\frac{(Cq^{1-n},dq^{n-1}/A,Cq^{1-n}/b,dq^n/AbZ)_{\infty}}
{(Z,d,q/b,Cd/AbZ)_{\infty}}
\prod_{i,j=1}^n
\frac{\left(\frac{x_i}{x_j}q,\frac{y_i}{y_j}a_jz_j\right)_\infty}
{\left(\frac{y_i}{y_j}q,\frac{x_i}{x_j}c_j\right)_\infty}\\\times
\sum_{k_1,\dots,k_n=-\infty}^{\infty}\Bigg(
\prod_{1\le i<j\le n}\left(\frac {y_iq^{k_i}-y_jq^{k_j}}{y_i-y_j}\right)
\prod_{i,j=1}^n\left(\frac{y_i}{y_j}a_jz_j\right)_{k_i}^{-1}
\prod_{i=1}^ny_i^{nk_i-|{\mathbf k}|}\\\times
\frac{(Aq^{1-n},AbZq^{1-n}/d)_{|{\mathbf k}|}}
{(Cq^{1-n})_{|{\mathbf k}|}}
(-1)^{(n-1)|{\mathbf k}|}
q^{-\binom{|{\mathbf k}|}2+n\sum_{i=1}^n\binom{k_i}2}
\left(\frac{dq^{n-1}}A\right)^{|\mathbf k|}\Bigg),
\end{multline}
provided that $|Cd/Ab|<|Z|<\left|q^{\frac
{n-1}2}x_j^{-n}\prod_{i=1}^nx_i\right|$ and
$|Cd/Ab|<|dq^{n-1}/A|<\left|q^{\frac
{n-1}2}y_j^{-n}\prod_{i=1}^ny_i\right|$
for $j=1,\dots,n$.
\end{Theorem}
\begin{proof}
We have, for
$\max(|Z|,|Cd/AbZ|,|dq^{n-1}/A|,|Cq^{1-n}/b|)<1$,
\begin{multline}\label{22t2gl}
{}_2\psi_2\!\left[\begin{matrix}Aq^{1-n},b\\
Cq^{1-n},d\end{matrix};q,Z\right]
=\frac{(AZq^{1-n},dq^{n-1}/A,Cq^{1-n}/b,dq^n/AbZ)_{\infty}}
{(Z,d,q/b,Cd/AbZ)_{\infty}}\\\times
{}_2\psi_2\!\left[\begin{matrix}Aq^{1-n},AbZq^{1-n}/d\\
AZq^{1-n},Cq^{1-n}\end{matrix};q,\frac{dq^{n-1}}A\right],
\end{multline}
by Bailey's $_2\psi_2$ transformation in \eqref{22tgl}.
Now we apply Lemma~\ref{lem1} to the $_2\psi_2$'s on the left and on the
right side of this transformation.
Specifically, we rewrite the $_2\psi_2$ on left side of \eqref{22t2gl}
by the $b_i\mapsto c_i$, $i=1,\dots,n$, and
\begin{equation*}
f(m)=\frac{(Aq^{1-n},b)_m}{(d)_m}Z^m
\end{equation*}
case of Lemma~\ref{lem1}.
The $_2\psi_2$ on the right side of \eqref{22t2gl} is rewritten
by the $b_i\mapsto a_iz_i$, $x_i\mapsto y_i$, $i=1,\dots,n$, and
\begin{equation*}
f(m)=\frac{(Aq^{1-n},AbZq^{1-n}/d)_m}{(Cq^{1-n})_m}
\left(\frac{dq^{n-1}}A\right)^m
\end{equation*}
case of Lemma~\ref{lem1}.
Finally, we divide both sides of the resulting equation by \eqref{divgl}
and simplify to obtain \eqref{an2psi22gl}.
\end{proof}

Next, we give two multivariable extensions of \eqref{22tgl}
which arise from Lemma~\ref{lem2}.

\begin{Theorem}[An $A_n$ $_2\psi_2$ transformation]\label{an2psi24}
Let $a_1,a_2,\dots,a_n$, $b$, $c_1,\dots,c_n$, $d$, $x_1,\dots,x_n$,
$y_1,\dots,y_n$, and $z$ be indeterminate, let $n\ge 1$, and suppose that
none of the denominators in \eqref{an2psi24gl} vanishes. Write
$A\equiv a_1\dots a_n$ and $C\equiv c_1\dots c_n$, for short. Then
\begin{multline}\label{an2psi24gl}
\sum_{k_1,\dots,k_n=-\infty}^{\infty}
\prod_{1\le i<j\le n}\left(\frac {x_iq^{k_i}-x_jq^{k_j}}{x_i-x_j}\right)
\prod_{i,j=1}^n\frac{\left(\frac{x_i}{x_j}a_j\right)_{k_i}}
{\left(\frac{x_i}{x_j}c_j\right)_{k_i}}\,
\frac{(b)_{|{\mathbf k}|}}{(d)_{|{\mathbf k}|}}
z^{|\mathbf k|}\\
=\frac{(Az,d/A,Cq^{1-n}/b,dq/Abz)_{\infty}}
{(z,d,q/b,Cdq^{1-n}/Abz)_{\infty}}
\prod_{i,j=1}^n
\frac{\left(\frac{y_i}{y_j}c_j,\frac{y_iq}{y_ja_i},
\frac{x_i}{x_j}q,\frac{x_ic_j}{x_ja_i}\right)_\infty}
{\left(\frac{x_i}{x_j}c_j,\frac{x_iq}{x_ja_i},
\frac{y_i}{y_j}q,\frac{y_ic_j}{y_ja_i}\right)_\infty}\\\times
\sum_{k_1,\dots,k_n=-\infty}^{\infty}
\prod_{1\le i<j\le n}\left(\frac {y_iq^{k_i}-y_jq^{k_j}}{y_i-y_j}\right)
\prod_{i,j=1}^n\frac{\left(\frac{y_i}{y_j}a_j\right)_{k_i}}
{\left(\frac{y_i}{y_j}c_j\right)_{k_i}}\,
\frac{(Abz/d)_{|{\mathbf k}|}}{(Az)_{|{\mathbf k}|}}
\left(\frac dA\right)^{|\mathbf k|},
\end{multline}
provided $|Cdq^{1-n}/Ab|<|z|<1$ and $|Cdq^{1-n}/Ab|<|d/A|<1$.
\end{Theorem}
\begin{proof}
We have, for
$\max(|z|,|Cdq^{1-n}/Abz|,|d/A|,|Cq^{1-n}/b|)<1$,
\begin{equation}\label{22t4gl}
{}_2\psi_2\!\left[\begin{matrix}A,b\\
Cq^{1-n},d\end{matrix};q,z\right]
=\frac{(Az,d/A,Cq^{1-n}/b,dq/Abz)_{\infty}}
{(z,d,q/b,Cdq^{1-n}/Abz)_{\infty}}\:
{}_2\psi_2\!\left[\begin{matrix}A,Abz/d\\
Az,Cq^{1-n}\end{matrix};q,\frac dA\right],
\end{equation}
by Bailey's $_2\psi_2$ transformation in \eqref{22tgl}.
Now we apply Lemma~\ref{lem2} to the $_2\psi_2$'s on the left and on the
right side of this transformation.
Specifically, we rewrite the $_2\psi_2$ on left side of \eqref{22t4gl}
by the $b_i\mapsto c_i$, $i=1,\dots,n$, and
\begin{equation*}
g(m)=\frac{(b)_m}{(d)_m}z^m
\end{equation*}
case of Lemma~\ref{lem2}.
The $_2\psi_2$ on the right side of \eqref{22t4gl} is rewritten
by the $b_i\mapsto c_i$, $x_i\mapsto y_i$, $i=1,\dots,n$, and
\begin{equation*}
g(m)=\frac{(Abz/d)_m}{(Az)_m}\left(\frac dA\right)^m
\end{equation*}
case of Lemma~\ref{lem2}.
Finally, we divide both sides of the resulting equation by
\begin{equation}\label{divgl2}
\frac{(q,Cq^{1-n}/A)_{\infty}}{(Cq^{1-n},q/A)_{\infty}}
\prod_{i,j=1}^n\frac{\left(\frac{x_i}{x_j}c_j,
\frac{x_iq}{x_ja_i}\right)_{\infty}}
{\left(\frac{x_i}{x_j}q,\frac{x_ic_j}{x_ja_i}\right)_{\infty}}
\end{equation}
and simplify to obtain \eqref{an2psi24gl}.
\end{proof}

\begin{Theorem}[An $A_n$ $_2\psi_2$ transformation]\label{an2psi25}
Let $a_1,\dots,a_n$, $b_1,\dots,b_n$, $c$, $d_1,\dots,d_n$, $x_1,\dots,x_n$,
$y_1,\dots,y_n$, and $z_1,\dots,z_n$ be indeterminate, let $n\ge 1$,
and suppose that none of the denominators in \eqref{an2psi25gl} vanishes.
Write $A\equiv a_1\dots a_n$, $B\equiv b_1\dots b_n$, $D\equiv d_1\dots d_n$,
and $Z\equiv z_1\dots z_n$, for short. Then
\begin{multline}\label{an2psi25gl}
\sum_{k_1,\dots,k_n=-\infty}^{\infty}
\prod_{1\le i<j\le n}\left(\frac {x_iq^{k_i}-x_jq^{k_j}}{x_i-x_j}\right)
\prod_{i,j=1}^n\frac{\left(\frac{x_i}{x_j}b_j\right)_{k_i}}
{\left(\frac{x_i}{x_j}d_j\right)_{k_i}}\,
\frac{(Aq^{1-n})_{|{\mathbf k}|}}{(c)_{|{\mathbf k}|}}
Z^{|\mathbf k|}\\
=\frac{(D/A,c/B)_{\infty}}{(Z,cD/ABZ)_{\infty}}
\prod_{i,j=1}^n
\frac{\left(\frac{y_i}{y_j}a_jz_j,\frac{y_id_iq}{y_ja_ib_iz_i},
\frac{x_i}{x_j}q,\frac{x_id_j}{x_jb_i}\right)_\infty}
{\left(\frac{x_i}{x_j}d_j,\frac{x_iq}{x_jb_i},
\frac{y_i}{y_j}q,\frac{y_id_ia_jz_j}{y_ja_ib_iz_i}\right)_\infty}\\\times
\sum_{k_1,\dots,k_n=-\infty}^{\infty}
\prod_{1\le i<j\le n}\left(\frac {y_iq^{k_i}-y_jq^{k_j}}{y_i-y_j}\right)
\prod_{i,j=1}^n\frac{\left(\frac{y_ia_jb_jz_j}{y_jd_j}\right)_{k_i}}
{\left(\frac{y_i}{y_j}a_jz_j\right)_{k_i}}\,
\frac{(Aq^{1-n})_{|{\mathbf k}|}}{(c)_{|{\mathbf k}|}}
\left(\frac DA\right)^{|\mathbf k|},
\end{multline}
provided $|cD/AB|<|Z|<1$ and $|cD/AB|<|D/A|<1$.
\end{Theorem}
\begin{proof}
We have, for
$\max(|Z|,|cD/ABZ|,|D/A|,|c/B|)<1$,
\begin{multline}\label{22t5gl}
{}_2\psi_2\!\left[\begin{matrix}Aq^{1-n},B\\
c,Dq^{1-n}\end{matrix};q,Z\right]\\
=\frac{(AZq^{1-n},D/A,c/B,Dq/ABZ)_{\infty}}
{(Z,Dq^{1-n},q/B,cD/ABZ)_{\infty}}\:
{}_2\psi_2\!\left[\begin{matrix}Aq^{1-n},ABZ/D\\
AZq^{1-n},c\end{matrix};q,\frac DA\right],
\end{multline}
by Bailey's $_2\psi_2$ transformation in \eqref{22tgl}.
Now we apply Lemma~\ref{lem2} to the $_2\psi_2$'s on the left and on the
right side of this transformation.
Specifically, we rewrite the $_2\psi_2$ on left side of \eqref{22t5gl}
by the $a_i\mapsto b_i$, $b_i\mapsto d_i$, $i=1,\dots,n$, and
\begin{equation*}
g(m)=\frac{(Aq^{1-n})_m}{(c)_m}Z^m
\end{equation*}
case of Lemma~\ref{lem2}.
The $_2\psi_2$ on the right side of \eqref{22t5gl} is rewritten
by the $a_i\mapsto a_ib_iz_i/d_i$, $b_i\mapsto a_iz_i$,
$x_i\mapsto y_i$, $i=1,\dots,n$, and
\begin{equation*}
g(m)=\frac{(Aq^{1-n})_m}{(c)_m}\left(\frac DA\right)^m
\end{equation*}
case of Lemma~\ref{lem2}.
Finally, we divide both sides of the resulting equation by
\begin{equation}\label{divgl1}
\frac{(q,Dq^{1-n}/B)_{\infty}}{(Dq^{1-n},q/B)_{\infty}}
\prod_{i,j=1}^n\frac{\left(\frac{x_i}{x_j}d_j,
\frac{x_iq}{x_jb_i}\right)_{\infty}}
{\left(\frac{x_i}{x_j}q,\frac{x_id_j}{x_jb_i}\right)_{\infty}}
\end{equation}
and simplify to obtain \eqref{an2psi25gl}.
\end{proof}

Finally, we provide two multivariable extensions of \eqref{22tgl1} which
arise from Lemmas~\ref{lem1} and \ref{lem2}, respectively.

\begin{Theorem}[An $A_n$ ${}_2\psi_2$ transformation]\label{an2psi23}
Let $a_1,a_2,\dots,a_n$, $b$, $c_1,\dots,c_n$, $d$, $x_1,\dots,x_n$,
$y_1,\dots,y_n$, and $z_1,\dots,z_n$ be indeterminate, let $n\ge 1$,
and suppose that none of the denominators in \eqref{an2psi23gl}
vanishes. Write $A\equiv a_1\dots a_n$, $C\equiv c_1\dots c_n$,
and $Z\equiv z_1\dots z_n$, for short. Then
\begin{multline}\label{an2psi23gl}
\sum_{k_1,\dots,k_n=-\infty}^{\infty}\Bigg(
\prod_{1\le i<j\le n}\left(\frac {x_iq^{k_i}-x_jq^{k_j}}{x_i-x_j}\right)
\prod_{i,j=1}^n\left(\frac{x_i}{x_j}c_j\right)_{k_i}^{-1}
\prod_{i=1}^nx_i^{nk_i-|{\mathbf k}|}\\\times
\frac{(Aq^{1-n},b)_{|{\mathbf k}|}}
{(d)_{|{\mathbf k}|}}(-1)^{(n-1)|{\mathbf k}|}
q^{-\binom{|{\mathbf k}|}2+n\sum_{i=1}^n\binom{k_i}2}
Z^{|\mathbf k|}\Bigg)\\
=\frac{(bZ,Cq/AbZ,dq^n/AbZ)_{\infty}}
{(q^n/A,q/b,d)_{\infty}}
\prod_{i,j=1}^n
\frac{\left(\frac{x_i}{x_j}q,\frac{y_i}{y_j}a_jz_j\right)_\infty}
{\left(\frac{y_i}{y_j}q,\frac{x_i}{x_j}c_j\right)_\infty}\\\times
\sum_{k_1,\dots,k_n=-\infty}^{\infty}\Bigg(
\prod_{1\le i<j\le n}\left(\frac {y_iq^{k_i}-y_jq^{k_j}}{y_i-y_j}\right)
\prod_{i,j=1}^n\left(\frac{y_i}{y_j}a_jz_j\right)_{k_i}^{-1}
\prod_{i=1}^ny_i^{nk_i-|{\mathbf k}|}\\\times
\frac{(AbZ/C,AbZq^{1-n}/d)_{|{\mathbf k}|}}
{(bZ)_{|{\mathbf k}|}}
(-1)^{(n-1)|{\mathbf k}|}
q^{-\binom{|{\mathbf k}|}2+n\sum_{i=1}^n\binom{k_i}2}
\left(\frac{Cd}{AbZ}\right)^{|\mathbf k|}\Bigg),
\end{multline}
provided that $|Cd/Ab|<|Z|<\left|q^{\frac
{n-1}2}x_j^{-n}\prod_{i=1}^nx_i\right|$ and
$|Cd/Ab|<|Cd/AbZ|<\left|q^{\frac
{n-1}2}y_j^{-n}\prod_{i=1}^ny_i\right|$
for $j=1,\dots,n$.
\end{Theorem}
\begin{proof}
We have, for
$\max(|Z|,|Cd/AbZ|)<1$,
\begin{multline}\label{22t3gl}
{}_2\psi_2\!\left[\begin{matrix}Aq^{1-n},b\\
Cq^{1-n},d\end{matrix};q,Z\right]
=\frac{(AZq^{1-n},bZ,Cq/AbZ,dq^n/AbZ)_{\infty}}
{(q^n/A,q/b,Cq^{1-n},d)_{\infty}}\\\times
{}_2\psi_2\!\left[\begin{matrix}AbZ/C,AbZq^{1-n}/d\\
AZq^{1-n},bZ\end{matrix};q,\frac{Cd}{AbZ}\right],
\end{multline}
by Bailey's $_2\psi_2$ transformation in \eqref{22tgl1}.
Now we apply Lemma~\ref{lem1} to the $_2\psi_2$'s on the left and on the
right side of this transformation.
Specifically, we rewrite the $_2\psi_2$ on left side of \eqref{22t3gl}
by the $b_i\mapsto c_i$, $i=1,\dots,n$, and
\begin{equation*}
f(m)=\frac{(Aq^{1-n},b)_m}{(d)_m}Z^m
\end{equation*}
case of Lemma~\ref{lem1}.
The $_2\psi_2$ on the right side of \eqref{22t3gl} is rewritten
by the $b_i\mapsto a_iz_i$, $x_i\mapsto y_i$, $i=1,\dots,n$, and
\begin{equation*}
f(m)=\frac{(AbZ/C,AbZq^{1-n}/d)_m}{(bZ)_m}
\left(\frac{Cd}{AbZ}\right)^m
\end{equation*}
case of Lemma~\ref{lem1}.
Finally, we divide both sides of the resulting equation by
\eqref{divgl} and simplify to obtain \eqref{an2psi23gl}.
\end{proof}

\begin{Theorem}[An $A_n$ $_2\psi_2$ transformation]\label{an2psi26}
Let $a_1,\dots,a_n$, $b_1,\dots,b_n$, $c$, $d_1,\dots,d_n$, $x_1,\dots,x_n$,
$y_1,\dots,y_n$, and $z_1,\dots,z_n$ be indeterminate, let $n\ge 1$,
and suppose that none of the denominators in \eqref{an2psi26gl} vanishes.
Write $A\equiv a_1\dots a_n$, $B\equiv b_1\dots b_n$, $D\equiv d_1\dots d_n$,
and $Z\equiv z_1\dots z_n$, for short. Then
\begin{multline}\label{an2psi26gl}
\sum_{k_1,\dots,k_n=-\infty}^{\infty}
\prod_{1\le i<j\le n}\left(\frac {x_iq^{k_i}-x_jq^{k_j}}{x_i-x_j}\right)
\prod_{i,j=1}^n\frac{\left(\frac{x_i}{x_j}b_j\right)_{k_i}}
{\left(\frac{x_i}{x_j}d_j\right)_{k_i}}\,
\frac{(Aq^{1-n})_{|{\mathbf k}|}}{(c)_{|{\mathbf k}|}}
Z^{|\mathbf k|}\\
=\frac{(BZ,cq^n/ABZ)_{\infty}}{(q^n/A,c)_{\infty}}
\prod_{i,j=1}^n
\frac{\left(\frac{y_i}{y_j}a_jz_j,\frac{y_id_iq}{y_ja_ib_iz_i},
\frac{x_i}{x_j}q,\frac{x_id_j}{x_jb_i}\right)_\infty}
{\left(\frac{x_i}{x_j}d_j,\frac{x_iq}{x_jb_i},
\frac{y_i}{y_j}q,\frac{y_id_ia_jz_j}{y_ja_ib_iz_i}\right)_\infty}\\\times
\sum_{k_1,\dots,k_n=-\infty}^{\infty}
\prod_{1\le i<j\le n}\left(\frac {y_iq^{k_i}-y_jq^{k_j}}{y_i-y_j}\right)
\prod_{i,j=1}^n\frac{\left(\frac{y_ia_jb_jz_j}{y_jd_j}\right)_{k_i}}
{\left(\frac{y_i}{y_j}a_jz_j\right)_{k_i}}\,
\frac{(ABZq^{1-n}/c)_{|{\mathbf k}|}}{(BZ)_{|{\mathbf k}|}}
\left(\frac{cD}{ABZ}\right)^{|\mathbf k|},
\end{multline}
provided $|cD/AB|<|Z|<1$.
\end{Theorem}
\begin{proof}
We have, for
$\max(|Z|,|cD/ABZ|)<1$,
\begin{multline}\label{22t6gl}
{}_2\psi_2\!\left[\begin{matrix}Aq^{1-n},B\\
c,Dq^{1-n}\end{matrix};q,Z\right]
=\frac{(AZq^{1-n},BZ,cq^n/ABZ,Dq/ABZ)_{\infty}}
{(q^n/A,q/B,c,Dq^{1-n})_{\infty}}\\\times
{}_2\psi_2\!\left[\begin{matrix}ABZq^{1-n}/c,ABZ/D\\
AZq^{1-n},BZ\end{matrix};q,\frac{cD}{ABZ}\right],
\end{multline}
by Bailey's $_2\psi_2$ transformation in \eqref{22tgl1}.
Now we apply Lemma~\ref{lem2} to the $_2\psi_2$'s on the left and on the
right side of this transformation.
Specifically, we rewrite the $_2\psi_2$ on left side of \eqref{22t6gl}
by the $a_i\mapsto b_i$, $b_i\mapsto d_i$, $i=1,\dots,n$, and
\begin{equation*}
g(m)=\frac{(Aq^{1-n})_m}{(c)_m}Z^m
\end{equation*}
case of Lemma~\ref{lem2}.
The $_2\psi_2$ on the right side of \eqref{22t6gl} is rewritten
by the $a_i\mapsto a_ib_iz_i/d_i$, $b_i\mapsto a_iz_i$,
$x_i\mapsto y_i$, $i=1,\dots,n$, and
\begin{equation*}
g(m)=\frac{(ABZq^{1-n}/c)_m}{(BZ)_m}\left(\frac{cD}{ABZ}\right)^m
\end{equation*}
case of Lemma~\ref{lem2}.
Finally, we divide both sides of the resulting equation by
\eqref{divgl1} and simplify to obtain \eqref{an2psi26gl}.
\end{proof}

\subsection{Some $\hbox{$\boldsymbol A_{\boldsymbol n}$}$
$\hbox{${}_{\boldsymbol 2}\boldsymbol\psi_{\boldsymbol 2}$}$ summations}
\label{subsec2}

Here, we work out (all) the $A_n$ extensions of the $_2\psi_2$ summation
in \eqref{22sgl} which arise from Lemmas~\ref{lem1} and \ref{lem2},
respectively.

First, we give two multivariable extensions of
\eqref{22sgl} which arise from Lemma~\ref{lem1}.

\begin{Theorem}[An $A_n$ $_2\psi_2$ summation]\label{an2psi2s}
Let $a$, $b$, $c_1,\dots,c_n$, and $x_1,\dots,x_n$ be indeterminate,
let $n\ge 1$, and suppose that
none of the denominators in \eqref{an2psi2sgl} vanishes. Then
\begin{multline}\label{an2psi2sgl}
\sum_{k_1,\dots,k_n=-\infty}^{\infty}\Bigg(
\prod_{1\le i<j\le n}\left(\frac {x_iq^{k_i}-x_jq^{k_j}}{x_i-x_j}\right)
\prod_{i,j=1}^n\left(\frac{x_i}{x_j}c_j\right)_{k_i}^{-1}
\prod_{i=1}^nx_i^{nk_i-|{\mathbf k}|}\\\times
\frac{(a,b)_{|{\mathbf k}|}}{(bq)_{|{\mathbf k}|}}(-1)^{(n-1)|{\mathbf k}|}
q^{-\binom{|{\mathbf k}|}2+n\sum_{i=1}^n\binom{k_i}2}
\left(\frac qa\right)^{|\mathbf k|}\Bigg)\\
=\frac{(q,bq/a,c_1\dots c_nq^{1-n}/b)_{\infty}}
{(q/a,bq,q/b)_{\infty}}
\prod_{i,j=1}^n
\frac{\left(\frac{x_i}{x_j}q\right)_\infty}
{\left(\frac{x_i}{x_j}c_j\right)_\infty},
\end{multline}
provided $|c_1\dots c_nq^{2-n}/a|<|q/a|<\left|q^{\frac
{n-1}2}x_j^{-n}\prod_{i=1}^nx_i\right|$
for $j=1,\dots,n$.
\end{Theorem}
\begin{proof}
We have, for
$\max(|q/a|,|c_1\dots c_nq^{1-n}|)<1$,
\begin{equation}\label{22s1gl}
{}_2\psi_2\!\left[\begin{matrix}a,b\\
c_1\dots c_nq^{1-n},bq\end{matrix};q,\frac qa\right]
=\frac{(q,q,bq/a,c_1\dots c_nq^{1-n}/b)_{\infty}}
{(q/a,bq,q/b,c_1\dots c_nq^{1-n})_{\infty}},
\end{equation}
by the $_2\psi_2$ summation in \eqref{22sgl}.
Now we apply Lemma~\ref{lem1} to the $_2\psi_2$ of this summation.
Specifically, we rewrite the $_2\psi_2$ in \eqref{22s1gl}
by the $b_i\mapsto c_i$, $i=1,\dots,n$, and
\begin{equation*}
f(m)=\frac{(a,b)_m}{(bq)_m}\left(\frac qa\right)^m
\end{equation*}
case of Lemma~\ref{lem1}.
Finally, we divide both sides of the resulting equation by
\eqref{divgl} and simplify to obtain \eqref{an2psi2sgl}.

For an alternative proof, set $z=q/a$ and $d=bq$ in Theorem~\ref{an2psi2}.
In this case the multilateral
series on the right side of \eqref{an2psi2gl} reduces to
\begin{multline*}
\sum_{\begin{smallmatrix}-\infty\le k_1,\dots,k_n\le\infty\\
|{\mathbf k}|=0\end{smallmatrix}}\Bigg(
\prod_{1\le i<j\le n}\left(\frac {y_iq^{k_i}-y_jq^{k_j}}{y_i-y_j}\right)
\prod_{i,j=1}^n\left(\frac{y_i}{y_j}c_j\right)_{k_i}^{-1}
\prod_{i=1}^ny_i^{nk_i-|{\mathbf k}|}\\\times
(-1)^{(n-1)|{\mathbf k}|}
q^{-\binom{|{\mathbf k}|}2+n\sum_{i=1}^n\binom{k_i}2}\Bigg)
=\frac{(c_1\dots c_nq^{1-n})_{\infty}}{(q)_{\infty}}
\prod_{i,j=1}^n\frac{\left(\frac{y_i}{y_j}q\right)_{\infty}}
{\left(\frac{y_i}{y_j}c_j\right)_{\infty}},
\end{multline*}
the last evaluation by the $m=0$ case of Proposition~\ref{an1psi1N}.
\end{proof}

\begin{Theorem}[An $A_n$ $_2\psi_2$ summation]\label{an2psi22s}
Let $a$, $b_1,\dots,b_n$, $c$, and $x_1,\dots,x_n$ be indeterminate,
let $n\ge 1$, and suppose that
none of the denominators in \eqref{an2psi22sgl} vanishes. Then
\begin{multline}\label{an2psi22sgl}
\sum_{k_1,\dots,k_n=-\infty}^{\infty}\Bigg(
\prod_{1\le i<j\le n}\left(\frac {x_iq^{k_i}-x_jq^{k_j}}{x_i-x_j}\right)
\prod_{i,j=1}^n\left(\frac{x_i}{x_j}b_jq\right)_{k_i}^{-1}
\prod_{i=1}^nx_i^{nk_i-|{\mathbf k}|}\\\times
\frac{(a,b_1\dots b_n)_{|{\mathbf k}|}}{(c)_{|{\mathbf k}|}}
(-1)^{(n-1)|{\mathbf k}|}
q^{-\binom{|{\mathbf k}|}2+n\sum_{i=1}^n\binom{k_i}2}
\left(\frac qa\right)^{|\mathbf k|}\Bigg)\\
=\frac{(q,b_1\dots b_nq/a,c/b_1\dots b_n)_{\infty}}
{(q/a,q/b_1\dots b_n,c)_{\infty}}
\prod_{i,j=1}^n
\frac{\left(\frac{x_i}{x_j}q\right)_\infty}
{\left(\frac{x_i}{x_j}b_jq\right)_\infty},
\end{multline}
provided $|cq/a|<|q/a|<\left|q^{\frac
{n-1}2}x_j^{-n}\prod_{i=1}^nx_i\right|$
for $j=1,\dots,n$.
\end{Theorem}
\begin{proof}
We have, for
$\max(|q/a|,|c|)<1$,
\begin{equation}\label{22s2gl}
{}_2\psi_2\!\left[\begin{matrix}a,b_1\dots b_n\\
c,b_1\dots b_nq\end{matrix};q,\frac qa\right]
=\frac{(q,q,b_1\dots b_nq/a,c/b_1\dots b_n)_{\infty}}
{(q/a,b_1\dots b_nq,q/b_1\dots b_n,c)_{\infty}},
\end{equation}
by the $_2\psi_2$ summation in \eqref{22sgl}.
Now we apply Lemma~\ref{lem1} to the $_2\psi_2$ of this summation.
Specifically, we rewrite the $_2\psi_2$ in \eqref{22s2gl}
by the $b_i\mapsto b_iq$, $i=1,\dots,n$, and
\begin{equation*}
f(m)=\frac{(a,b_1\dots b_n)_m}{(c)_m}\left(\frac qa\right)^m
\end{equation*}
case of Lemma~\ref{lem1}.
Finally, we divide both sides of the resulting equation by
\begin{equation*}
\frac{(q)_{\infty}}{(b_1\dots b_nq^{1-n})_{\infty}}
\prod_{i,j=1}^n\frac{\left(\frac{x_i}{x_j}b_j\right)_{\infty}}
{\left(\frac{x_i}{x_j}q\right)_{\infty}}
\end{equation*}
and simplify to obtain \eqref{an2psi22sgl}.

For an alternative proof, set $c_i=b_iq$, $z_i=q/a_i$, $i=1,\dots,n$, and
$b\mapsto b_1\dots b_n$ in Theorem~\ref{an2psi23}.
In this case the multilateral
series on the right side of \eqref{an2psi23gl} is terminated from
below and from above and reduces just to one term, 1.
In the resulting equation, we replace
$A$ by $aq^{n-1}$ and $d$ by $c$.
\end{proof}

Finally, we give four multivariable extensions of
\eqref{22sgl} which arise from Lemma~\ref{lem2}.

\begin{Theorem}[An $A_n$ $_2\psi_2$ summation]\label{an2psi23s}
Let $a_1,\dots,a_n$, $b$, $c_1,\dots,c_n$, and $x_1,\dots,x_n$ be
indeterminate, let $n\ge 1$, and suppose that
none of the denominators in \eqref{an2psi23sgl} vanishes. Then
\begin{multline}\label{an2psi23sgl}
\sum_{k_1,\dots,k_n=-\infty}^{\infty}
\prod_{1\le i<j\le n}\left(\frac {x_iq^{k_i}-x_jq^{k_j}}{x_i-x_j}\right)
\prod_{i,j=1}^n\frac{\left(\frac{x_i}{x_j}a_j\right)_{k_i}}
{\left(\frac{x_i}{x_j}c_j\right)_{k_i}}\,
\frac{(b)_{|{\mathbf k}|}}{(bq)_{|{\mathbf k}|}}
\left(\frac q{a_1\dots a_n}\right)^{|\mathbf k|}\\
=\frac{(q,bq/a_1\dots a_n,c_1\dots c_nq^{1-n}/b)_{\infty}}
{(bq,q/b,c_1\dots c_nq^{1-n}/a_1\dots a_n)_{\infty}}
\prod_{i,j=1}^n
\frac{\left(\frac{x_i}{x_j}q,\frac{x_ic_j}{x_ja_i}\right)_\infty}
{\left(\frac{x_i}{x_j}c_j,\frac{x_iq}{x_ja_i}\right)_\infty},
\end{multline}
provided $\max(|c_1\dots c_nq^{1-n}|,|q/a_1\dots a_n|)<1$.
\end{Theorem}
\begin{proof}
We have, for
$\max(|q/a_1\dots a_n|,|c_1\dots c_nq^{1-n}|)<1$,
\begin{equation}\label{22s3gl}
{}_2\psi_2\!\left[\begin{matrix}a_1\dots a_n,b\\
c_1\dots c_nq^{1-n},bq\end{matrix};q,\frac q{a_1\dots a_n}\right]
=\frac{(q,q,bq/a_1\dots a_n,c_1\dots c_nq^{1-n}/b)_{\infty}}
{(q/a_1\dots a_n,bq,q/b,c_1\dots c_nq^{1-n})_{\infty}},
\end{equation}
by the $_2\psi_2$ summation in \eqref{22sgl}.
Now we apply Lemma~\ref{lem2} to the $_2\psi_2$ of this summation.
Specifically, we rewrite the $_2\psi_2$ in \eqref{22s3gl}
by the $b_i\mapsto c_i$, $i=1,\dots,n$, and
\begin{equation*}
g(m)=\frac{(b)_m}{(bq)_m}\left(\frac q{a_1\dots a_n}\right)^m
\end{equation*}
case of Lemma~\ref{lem2}.
Finally, we divide both sides of the resulting equation by
\eqref{divgl2} and simplify to obtain \eqref{an2psi23sgl}.

For an alternative proof, set $z=q/a_1\dots a_n$ and $d=bq$ in
Theorem~\ref{an2psi24}.
In this case the multilateral
series on the right side of \eqref{an2psi24gl} reduces to
\begin{multline*}
\sum_{\begin{smallmatrix}-\infty\le k_1,\dots,k_n\le\infty\\
|{\mathbf k}|=0\end{smallmatrix}}
\prod_{1\le i<j\le n}\left(\frac {y_iq^{k_i}-y_jq^{k_j}}{y_i-y_j}\right)
\prod_{i,j=1}^n\frac{\left(\frac{y_i}{y_j}a_j\right)_{k_i}}
{\left(\frac{y_i}{y_j}c_j\right)_{k_i}}\\
=\frac{(c_1\dots c_nq^{1-n},q/a_1\dots a_n)_{\infty}}
{(q,c_1\dots c_nq^{1-n}/a_1\dots a_n)_{\infty}}
\prod_{i,j=1}^n\frac{\left(\frac{y_i}{y_j}q,
\frac{y_ic_j}{y_ja_i}\right)_{\infty}}
{\left(\frac{y_i}{y_j}c_j,
\frac{y_iq}{y_ja_i}\right)_{\infty}},
\end{multline*}
the last evaluation by Theorem~\ref{an1psi1cN}.
\end{proof}

\begin{Theorem}[An $A_n$ $_2\psi_2$ summation]\label{an2psi24s}
Let $a$, $b_1,\dots,b_n$, $c$, and $x_1,\dots,x_n$ be
indeterminate, let $n\ge 1$, and suppose that
none of the denominators in \eqref{an2psi24sgl} vanishes. Then
\begin{multline}\label{an2psi24sgl}
\sum_{k_1,\dots,k_n=-\infty}^{\infty}
\prod_{1\le i<j\le n}\left(\frac {x_iq^{k_i}-x_jq^{k_j}}{x_i-x_j}\right)
\prod_{i,j=1}^n\frac{\left(\frac{x_i}{x_j}b_j\right)_{k_i}}
{\left(\frac{x_i}{x_j}b_jq\right)_{k_i}}\,
\frac{(a)_{|{\mathbf k}|}}{(c)_{|{\mathbf k}|}}
\left(\frac qa\right)^{|\mathbf k|}\\
=\frac{(b_1\dots b_nq/a,c/b_1\dots b_n)_{\infty}}
{(q/a,c)_{\infty}}
\prod_{i,j=1}^n
\frac{\left(\frac{x_i}{x_j}q,\frac{x_ib_j}{x_jb_i}q\right)_\infty}
{\left(\frac{x_i}{x_j}b_jq,\frac{x_iq}{x_jb_i}\right)_\infty},
\end{multline}
provided $\max(|c|,|q/a|)<1$.
\end{Theorem}
\begin{proof}
We utilize the $_2\psi_2$ summation in \eqref{22s1gl} and
apply Lemma~\ref{lem2} to the $_2\psi_2$ in that summation.
Specifically, we rewrite the $_2\psi_2$ in \eqref{22s1gl}
by the $a_i\mapsto b_i$, $b_i\mapsto b_iq$, $i=1,\dots,n$, and
\begin{equation*}
g(m)=\frac{(a)_m}{(c)_m}\left(\frac qa\right)^m
\end{equation*}
case of Lemma~\ref{lem2}.
Finally, we divide both sides of the resulting equation by
\begin{equation*}
\frac{(q,q)_{\infty}}{(b_1\dots b_nq,q/b_1\dots b_n)_{\infty}}
\prod_{i,j=1}^n\frac{\left(\frac{x_i}{x_j}b_jq,
\frac{x_iq}{x_jb_i}\right)_{\infty}}
{\left(\frac{x_i}{x_j}q,\frac{x_ib_j}{x_jb_i}q\right)_{\infty}}
\end{equation*}
and simplify to obtain \eqref{an2psi24sgl}.

For an alternative proof, set $z_i=q/a_i$, and
$d_i=b_iq$, $i=1,\dots,n$, in Theorem~\ref{an2psi25}.
In this case the multilateral
series on the right side of \eqref{an2psi25gl} is terminated from
below and from above and reduces just to
one term, 1. In the resulting summation, replace $A$ by $aq^{n-1}$.
\end{proof}

\begin{Theorem}[An $A_n$ $_2\psi_2$ summation]\label{an2psi25s}
Let $a$, $b_1,\dots,b_n$, $c_1,\dots,c_n$, and $x_1,\dots,x_n$ be
indeterminate, let $n\ge 1$, and suppose that
none of the denominators in \eqref{an2psi25sgl} vanishes. Then
\begin{multline}\label{an2psi25sgl}
\sum_{k_1,\dots,k_n=-\infty}^{\infty}
\prod_{1\le i<j\le n}\left(\frac {x_iq^{k_i}-x_jq^{k_j}}{x_i-x_j}\right)
\prod_{i,j=1}^n\frac{\left(\frac{x_i}{x_j}b_j\right)_{k_i}}
{\left(\frac{x_i}{x_j}c_j\right)_{k_i}}\,
\frac{(a)_{|{\mathbf k}|}}{(b_1\dots b_nq)_{|{\mathbf k}|}}
\left(\frac qa\right)^{|\mathbf k|}\\
=\frac{(q,b_1\dots b_nq/a)_{\infty}}
{(q/a,b_1\dots b_nq)_{\infty}}
\prod_{i,j=1}^n
\frac{\left(\frac{x_i}{x_j}q,\frac{x_ic_j}{x_jb_i}\right)_\infty}
{\left(\frac{x_i}{x_j}c_j,\frac{x_iq}{x_jb_i}\right)_\infty},
\end{multline}
provided $\max(|c_1\dots c_nq^{1-n}|,|q/a|)<1$.
\end{Theorem}
\begin{proof}
Write $B\equiv b_1\dots b_n$ and $C\equiv c_1\dots c_n$. We have, for
$\max(|q/a|,|Cq^{1-n}|)<1$,
\begin{equation}\label{22s5gl}
{}_2\psi_2\!\left[\begin{matrix}a,B\\
Cq^{1-n},Bq\end{matrix};q,\frac qa\right]
=\frac{(q,q,Bq/a,Cq^{1-n}/B)_{\infty}}
{(q/a,Bq,q/B,Cq^{1-n})_{\infty}},
\end{equation}
by the $_2\psi_2$ summation in \eqref{22sgl}.
Now we apply Lemma~\ref{lem2} to the $_2\psi_2$ of this summation.
Specifically, we rewrite the $_2\psi_2$ in \eqref{22s5gl}
by the $a_i\mapsto b_i$, $b_i\mapsto c_i$, $i=1,\dots,n$, and
\begin{equation*}
g(m)=\frac{(a)_m}{(Bq)_m}\left(\frac qa\right)^m
\end{equation*}
case of Lemma~\ref{lem2}.
Finally, we divide both sides of the resulting equation by
\begin{equation*}
\frac{(q,Cq^{1-n}/B)_{\infty}}
{(Cq^{1-n},q/B)_{\infty}}
\prod_{i,j=1}^n\frac{\left(\frac{x_i}{x_j}c_j,
\frac{x_iq}{x_jb_i}\right)_{\infty}}
{\left(\frac{x_i}{x_j}q,\frac{x_ic_j}{x_jb_i}\right)_{\infty}}
\end{equation*}
and simplify to obtain \eqref{an2psi25sgl}.

For an alternative proof, set $z_i=q/a_i$, $i=1,\dots,n$,
and $c=b_1\dots b_nq$ in Theorem~\ref{an2psi26}.
In this case the multilateral
series on the right side of \eqref{an2psi26gl} is terminated from
below and from above and reduces just to
one term, 1. In the resulting summation, replace $d_i$ by $c_i$,
$i=1,\dots,n$, and $A$ by $aq^{n-1}$.
\end{proof}

\begin{Theorem}[An $A_n$ $_2\psi_2$ summation]\label{an2psi26s}
Let $a_1,\dots,a_n$, $b_1,\dots,b_n$, $c$, and $x_1,\dots,x_n$ be
indeterminate, let $n\ge 1$, and suppose that
none of the denominators in \eqref{an2psi26sgl} vanishes. Then
\begin{multline}\label{an2psi26sgl}
\sum_{k_1,\dots,k_n=-\infty}^{\infty}
\prod_{1\le i<j\le n}\left(\frac {x_iq^{k_i}-x_jq^{k_j}}{x_i-x_j}\right)
\prod_{i,j=1}^n\frac{\left(\frac{x_i}{x_j}a_j\right)_{k_i}}
{\left(\frac{x_i}{x_j}b_jq\right)_{k_i}}\,
\frac{(b_1\dots b_n)_{|{\mathbf k}|}}{(c)_{|{\mathbf k}|}}
\left(\frac q{a_1\dots a_n}\right)^{|\mathbf k|}\\
=\frac{(q,c/b_1\dots b_n)_{\infty}}
{(q/b_1\dots b_n,c)_{\infty}}
\prod_{i,j=1}^n
\frac{\left(\frac{x_i}{x_j}q,\frac{x_ib_jq}{x_ja_i}\right)_\infty}
{\left(\frac{x_i}{x_j}b_jq,\frac{x_iq}{x_ja_i}\right)_\infty},
\end{multline}
provided $\max(|c|,|q/a_1\dots a_n|)<1$.
\end{Theorem}
\begin{proof}
Write $A\equiv a_1\dots a_n$ $B\equiv b_1\dots b_n$. We have, for
$\max(|q/A|,|c|)<1$,
\begin{equation}\label{22s6gl}
{}_2\psi_2\!\left[\begin{matrix}A,B\\
c,Bq\end{matrix};q,\frac qA\right]
=\frac{(q,q,Bq/A,c/B)_{\infty}}
{(q/A,Bq,q/B,c)_{\infty}},
\end{equation}
by the $_2\psi_2$ summation in \eqref{22sgl}.
Now we apply Lemma~\ref{lem2} to the $_2\psi_2$ of this summation.
Specifically, we rewrite the $_2\psi_2$ in \eqref{22s6gl}
by the $b_i\mapsto b_iq$, $i=1,\dots,n$, and
\begin{equation*}
g(m)=\frac{(B)_m}{(c)_m}\left(\frac qA\right)^m
\end{equation*}
case of Lemma~\ref{lem2}.
Finally, we divide both sides of the resulting equation by
\begin{equation*}
\frac{(q,Bq/A)_{\infty}}
{(Bq,q/A)_{\infty}}
\prod_{i,j=1}^n\frac{\left(\frac{x_i}{x_j}b_jq,
\frac{x_iq}{x_ja_i}\right)_{\infty}}
{\left(\frac{x_i}{x_j}q,\frac{x_ib_jq}{x_ja_i}\right)_{\infty}}
\end{equation*}
and simplify to obtain \eqref{an2psi26sgl}.

For an alternative proof, set $z_i=q^{\frac 1n}/b_i$,
and $d_i=a_iq^{\frac 1n}$, $i=1,\dots,n$, in Theorem~\ref{an2psi26}.
In this case the multilateral
series on the right side of \eqref{an2psi26gl} is terminated from
below and from above and reduces just to
one term, 1. In the resulting summation, replace $a_i$ by
$b_iq^{1-\frac 1n}$, and $b_i$ by $a_i$, for $i=1,\dots,n$.
\end{proof}

\end{document}